\documentclass[10pt,compsocconf]{IEEEtran}

\usepackage[cmex10]{amsmath}
\usepackage{subfig}
\usepackage{url}

\usepackage{amssymb,bm,graphicx,cite}
\usepackage[mathscr]{eucal}
  \renewcommand{\b}{\mathbf}
  \newcommand{\R}{\mathbb{R}}
  \newcommand{\C}{\bm{\mathscr{C}}}
  \newcommand{\G}{\mathscr{G}}
  \newcommand{\D}{\bm{\Delta}}
    \newcommand{\dx}{{d_x}}
    \newcommand{\da}{{d_{\alpha}}}
  \newcommand{\CA}{\bm{\mathscr{A}}}
  \newcommand{\CD}{\bm{\mathscr{D}}}

\DeclareMathOperator*{\argmin}{arg\,min}
\begin{document}
\title{Collaborative Filtering via Group-Structured Dictionary Learning}
\author{\IEEEauthorblockN{Zolt{\'a}n Szab{\'o}\IEEEauthorrefmark{1},
Barnab{\'a}s P{\'o}czos\IEEEauthorrefmark{2},
and Andr{\'a}s L{\H{o}}rincz\IEEEauthorrefmark{1}}\\
\IEEEauthorblockA{\IEEEauthorrefmark{1}Faculty of Informatics, E{\"o}tv{\"o}s Lor\'{a}nd University, P{\'a}zm{\'a}ny P{\'e}ter s{\'e}t{\'a}ny 1/C, H-1117 Budapest, Hungary\\
Email: \url{szzoli@cs.elte.hu}, \url{andras.lorincz@elte.hu}, Web: \url{http://nipg.inf.elte.hu}}\\
\IEEEauthorblockA{\IEEEauthorrefmark{2}Carnegie Mellon University, Robotics Institute, 5000 Forbes Ave, Pittsburgh, PA 15213\\
Email: \url{bapoczos@cs.cmu.edu}, Web: \url{http://www.autonlab.org}}%
\thanks{A compressed version of the paper has been accepted for publication at the $10^{th}$ International Conference on Latent Variable Analysis and Source Separation (LVA/ICA 2012).}}
\maketitle

\begin{abstract}
Structured sparse coding and the related structured dictionary learning problems are novel research areas in machine learning.
In this paper we present a new application of structured dictionary learning for collaborative filtering based
recommender systems. Our extensive numerical experiments demonstrate that the presented technique outperforms its
state-of-the-art competitors and has several advantages over approaches that do not put structured constraints on the
dictionary elements.
\end{abstract}


\section{Introduction}\label{sec:introduction}

The proliferation of online services and the thriving electronic commerce overwhelms us with
alternatives in our daily lives. To handle this information overload and to help users in efficient decision making, recommender systems (RS) have been designed. The goal of RSs is to recommend personalized items for online users when they need to choose among several items. Typical problems include recommendations for which movie to watch,
which jokes/books/news to read, which hotel to stay at, or which songs to listen to.

One of the most popular approaches in the field of recommender systems is \emph{collaborative filtering} (CF). The underlying idea of CF is very simple: Users generally express their tastes in an explicit way by rating the items. CF tries to estimate the users' preferences based on the ratings they have already made on items
and based on the ratings of other, similar users.
For a recent review on recommender systems and collaborative filtering, see e.g., \cite{ricci11recommender}.

Novel advances on CF show that \emph{dictionary learning} based approaches can be efficient for making predictions about users' preferences \cite{takacs09scalable}. The dictionary learning based approach assumes that (i) there is a latent, unstructured feature space (hidden representation) behind the users' ratings, and (ii) a rating of an item is equal to the product of the item and the user's feature. To increase the generalization capability, usually $\ell_2$ regularization is introduced both for the dictionary and for the users' representation.

There are several problems that belong to the task of dictionary learning \cite{yaghoobi09dictionary}, a.k.a. matrix factorization \cite{witten09penalized}. This set of problems includes, for example, (sparse) principal component analysis \cite{zou06sparse}, independent component analysis
\cite{hyvarinen01independent}, independent subspace analysis \cite{cardoso98multidimensional}, non-negative matrix factorization \cite{lee00algorithms}, and
\emph{structured dictionary} learning, which will be the target of our paper.

One predecessor of the structured dictionary learning problem is the \emph{sparse coding} task \cite{tropp10computational}, which is a considerably simpler problem. Here the dictionary is already given, and we assume that the observations can be approximated well enough using only a few dictionary elements. Although finding the solution
that uses the minimal number of dictionary elements is NP hard in general \cite{natarajan95sparse}, there exist efficient approximations.
One prominent example is the Lasso approach \cite{tibshirani96regression}, which applies convex $\ell_1$ relaxation to the
code words. Lasso does not enforce any \emph{group} structure on the components of the representation (covariates).

However, using \emph{structured sparsity}, that is, forcing
different kind of structures (e.g., disjunct groups, trees) on the sparse codes can
lead to increased performances in several applications. Indeed, as it has been theoretically proved recently
structured sparsity can ease feature selection \cite{huang10benefit,yuan06model}, and
makes possible robust compressed sensing with substantially decreased observation number \cite{baraniuk10model}.
Many other real life applications also confirm the benefits of structured sparsity, for example
(i) automatic image annotation \cite{zhang10automatic}, (ii) group-structured feature selection for micro array data processing
\cite{zhao09composite,jacob09group,kim10tree,rapaport08classification},
(iii) multi-task learning problems (a.k.a. transfer learning) \cite{obozinski10joint,kim10scalable,rakotomamonjy11review},
(iv) multiple kernel learning \cite{szafranski10composite,aflalo11variable},
(v) face recognition \cite{elhamifar11robust}, and (vi) structure learning in graphical models \cite{schmidt10convex,jalali11learning}.
For an excellent review on structured sparsity, see \cite{bach11optimization}.

All the above mentioned examples only consider the structured sparse coding problem, where we assume that the dictionary is already given and available to us. A more interesting (and challenging) problem is the combination of these two tasks, i.e., learning the best structured dictionary and structured representation.
This is the \emph{structured dictionary learning} (SDL) problem. SDL is more difficult; one can find only few solutions in the literature
\cite{jenatton11proximal,jenatton10structured,mairal11convex,rosenblum10dictionary,koray09learning,silva11blind}.
This novel field is appealing for (i) transformation invariant feature extraction \cite{koray09learning},
(ii) image denoising/inpainting \cite{jenatton11proximal,mairal11convex,silva11blind}, (iii) background subtraction \cite{mairal11convex},
(iv) analysis of text corpora \cite{jenatton11proximal}, and (v) face recognition \cite{jenatton10structured}.

\textbf{Our goal} is to extend the application domain of SDL in the direction of collaborative filtering. With respect to CF,
further constraints appear for SDL since (i) online learning is desired and (ii) missing information is typical.
There are good reasons for them: novel items/users may appear and user preferences may change over time.
Adaptation to users also motivate online methods. Online methods have the additional advantage with respect to offline ones
that they can process
more instances in the same amount of time, and in many cases this can lead to increased performance.
For a theoretical proof of this claim, see \cite{bottou05on-line}. Furthermore, users can evaluate only a small portion of the
available items, which leads to incomplete observations, missing rating values. 
In order to cope with these constraints of the collaborative filtering problem, we will use a 
novel extension of the structured dictionary learning problem, the so-called online group-structured dictionary learning (OSDL)  \cite{szabo11online}. OSDL allows (i)
overlapping group structures with (ii)
non-convex sparsity inducing regularization, (iii) partial observation (iv) in an online framework.

Our paper is structured as follows: We briefly review the OSDL problem, its cost function, and optimization method in Section~\ref{sec:OSDL problem}. We cast the CF problem as an OSDL task in Section~\ref{sec:OSDL via CF}. Numerical results are presented in Section~\ref{sec:numerical results}. Conclusions are drawn in Section~\ref{sec:conclusions}.

\textbf{Notations.} Vectors ($\b{a}$) and matrices ($\b{A}$) are denoted by bold letters. $diag(\b{a})$ represents the diagonal
matrix with coordinates of vector $\b{a}$ in its diagonal. The $i^{th}$ coordinate of vector $\b{a}$ is $a_i$. Notation
$|\cdot|$ means the number of elements of a set and the absolute value for a real number. For set $O\subseteq\{1,\ldots,d\}$,
$\b{a}_O\in\R^{|O|}$ denotes the coordinates of vector $\b{a}\in\R^d$ in $O$. For matrix
$\b{A}\in\R^{d\times D}$, $\b{A}_O\in\R^{|O|\times D}$
stands for the restriction of matrix $\b{A}$ to the rows $O$. $\b{I}$ and  $\b{0}$ denote the identity and the null matrices, respectively. $\b{A}^T$ is the
transposed form of $\b{A}$. For a vector, the $\max$ operator acts coordinate-wise. The  $\ell_{p}$ (quasi-)norm of vector $\b{a}\in\R^{d}$ is $\|\b{a}\|_{p}=(\sum_{i=1}^d|a_i|^{p})^{\frac{1}{p}}$ ($p>0$). $S_p^d=\{\b{a}\in \R^d:\|\b{a}\|_p\le 1\}$ denotes the $\ell_p$ unit sphere in $\R^d$. The point-wise and scalar products of $\b{a},\b{b}\in\R^d$ are denoted by $\b{a}\circ\b{b}=[a_1b_1;\ldots;a_db_d]$ and by $\left<\b{a},\b{b}\right> = \b{a}^T\b{b}$, respectively.
For a set system $\G$, the coordinates of vector $\b{a}\in\R^{|\G|}$ are denoted by $a^G$ ($G\in\G$), that is, $\b{a}=(a^G)_{G\in\G}$. $\Pi_{\C}(\b{x})=\mathrm{argmin}_{\b{c}\in\C}\|\b{x}-\b{c}\|_2$ is the projection of point $\b{x}\in\R^d$ to the convex closed set $\C\subseteq\R^d$. Partial derivative of function $h$ w.r.t. variable $\b{x}$ in $\b{x}_0$ is $\frac{\partial h}{\partial \b{x}}(\b{x}_0)$. The non-negative ortant of $\R^d$ is $\R^{d}_{+}=\{\b{x}\in\R^d:x_i\ge 0\, (\forall i)\}$. For sets, $\times$ and $\backslash$ denote direct product and difference, respectively.

\section{The OSDL Problem}\label{sec:OSDL problem}
In this section we briefly review the OSDL approach, which will be our major tool to solve the CF problem. The OSDL cost function is treated in Section~\ref{sec:OSDL:cost},
its optimization idea is detailed in Section~\ref{sec:OSDL:optimization}.

\subsection{Cost Function}\label{sec:OSDL:cost}
The online group-structured dictionary learning (OSDL) task is defined with the following quantities. Let the dimension of the observations
be denoted by $\dx$. Assume that in each time instant ($i=1,2,\ldots$) a set $O_i\subseteq\{1,\ldots,\dx\}$ is given,
that is, we know which coordinates are observable at time $i$, and the observation is $\b{x}_{O_i}$. Our goal is to find a
dictionary $\b{D}\in \R^{\dx\times \da}$ that can approximate the observations $\b{x}_{O_i}$ well from the linear combination of
its columns. The columns of $\b{D}$ are assumed to belong to a closed, convex, and bounded set $\CD=\times_{i=1}^{\da}\CD_i$.
To formulate the cost of dictionary $\b{D}$, first a \emph{fixed} time instant $i$, observation $\b{x}_{O_i}$, dictionary $\b{D}$ is considered, and
the hidden representation $\bm{\alpha}_i$ associated to this $(\b{x}_{O_i},\b{D},O_i)$ triple is defined. Representation $\bm{\alpha}_i$ is allowed to belong to
a closed, convex set $\CA\subseteq \R^{\da}$ ($\bm{\alpha}_i\in\CA$) with certain structural constraints.
The structural constraint on $\bm{\alpha}_i$ are expressed by making use of a given $\G$ group structure, which is a set system (also called hypergraph) on
$\{1,\ldots,\da\}$. It is also assumed that weight vectors $\b{d}^{G}\in\R^{\da}$ ($G\in\G$)  are available for us and that they are
positive on $G$ and $0$ otherwise. Representation $\bm{\alpha}$ belonging to a triple $(\b{x}_O,\b{D},O)$ is defined as the solution
of the structured sparse coding task
\begin{align}
l(\b{x}_O,\b{D}_O) &= l_{\CA,\kappa,\G,\left\{\b{d}^{G}\right\}_{G\in\G},\eta}(\b{x}_O,\b{D}_O)\\
                   &= \min_{\bm{\alpha}\in\CA}\left[\frac{1}{2}\left\|\b{x}_O-\b{D}_O\bm{\alpha}\right\|^2_2+\kappa\Omega(\bm{\alpha})\right],\label{eq:l}
\end{align}
where $l(\b{x}_O,\b{D}_O)$ denotes the loss, $\kappa>0$,  and
\begin{equation}
\Omega(\b{y}) = \Omega_{\G,\left\{\b{d}^{G}\right\}_{G\in\G},\eta}(\b{y}) = \|(\|\b{d}^{G}\circ \b{y}\|_2)_{G\in\G}\|_{\eta}\label{eq:Omega}
\end{equation}
is the structured regularizer associated to $\G$ and $\{\b{d}^{G}\}_{G\in\G}$, $\eta\in(0,2)$. Here, the first term of \eqref{eq:l} is responsible
for the quality of approximation on the observed coordinates, whereas for $\eta\le 1$ the other term [\eqref{eq:Omega}] constrains the solution
according to the group structure $\G$ similarly to the sparsity inducing regularizer $\Omega$ in \cite{jenatton10structured}: it
eliminates the terms $\|\b{d}^{G}\circ\b{y}\|_2$ $(G\in\G)$ by means of $\left\|\cdot\right\|_{\eta}$. The OSDL problem is defined as
the minimization of the cost function:
\begin{equation}
      \min_{\b{D}\in\CD}f_t(\b{D}):=\frac{1}{\sum_{j=1}^t(j/t)^{\rho}}\sum_{i=1}^t\left(\frac{i}{t}\right)^{\rho}l(\b{x}_{O_i},\b{D}_{O_i}),\label{eq:f_t(D)-2}
\end{equation}
that is, the goal is to minimize the average loss belonging to the dictionary, where $\rho$ is a non-negative forgetting factor.
If $\rho=0$, the classical average $f_t(\b{D})=\frac{1}{t}\sum_{i=1}^tl(\b{x}_{O_i},\b{D}_{O_i})$ is recovered.

As an example, let $\CD_i=S^{\dx}_2$ ($\forall i$), $\CA=\R^{\da}$. In this case, columns of $\b{D}$ are restricted to the Euclidean
unit sphere and we have no constraints for $\bm{\alpha}$. Now, let $|\G|=\da$ and $\G=\{desc_1,\ldots,desc_{\da}\}$,
where $desc_i$ represents the $i^{th}$ node and its children in a fixed tree. Then the
 coordinates $\alpha_i$ are searched in a hierarchical tree structure and the hierarchical dictionary  $\b{D}$ is optimized accordingly.

\subsection{Optimization}\label{sec:OSDL:optimization}
Optimization of cost function \eqref{eq:f_t(D)-2} is equivalent to the joint optimization of dictionary $\b{D}$ and representation $\{\bm{\alpha}_i\}_{i=1}^t$:
\begin{equation}
\argmin\limits_{\b{D}\in\CD, \{\bm{\alpha}_i\in\CA\}_{i=1}^t}\hspace*{-.3cm}f_t(\b{D},\{\bm{\alpha}_i\}_{i=1}^t),\label{eq:f(D,{alfa}):min-only}
\end{equation}
where
{\small\begin{equation}
  f_t=\frac{1}{\sum_{j=1}^t(j/t)^{\rho}}\sum_{i=1}^t\left(\frac{i}{t}\right)^{\rho} \left[\frac{1}{2}\left\|\b{x}_{O_i}-\b{D}_{O_i}\bm{\alpha}_i\right\|^2_2+\kappa\Omega(\bm{\alpha}_i)\right]. \label{eq:f(D,{alfa})}
\end{equation}}
$\b{D}$ is optimized by using the sequential observations $\b{x}_{O_i}$ online in an alternating manner:
\begin{enumerate}
  \item  The actual dictionary estimation $\b{D}_{t-1}$ and sample $\b{x}_{O_t}$ is used to optimize \eqref{eq:l} for  representation $\bm{\alpha}_t$.
   \item For the estimated representations $\{\bm{\alpha}_i\}_{i=1}^t$, the dictionary estimation $\b{D}_{t}$ is derived from the quadratic optimization problem
    \begin{equation}
      \hat{f}_t(\b{D}_t) = \min_{\b{D}\in \CD} f_t(\b{D},\{\bm{\alpha}_i\}_{i=1}^t).\label{eq:task2(D-opt)}
    \end{equation}
\end{enumerate}

\subsubsection{Representation optimization ($\bm{\alpha}$)}
Note that \eqref{eq:l} is a non-convex optimization problem with respect to $\bm{\alpha}$. The variational properties
of norm $\eta$ can be used to overcome this problem. One can show, alike to \cite{jenatton10structured}, that by introducing an
auxiliary variable $\bm{z}\in\R_+^{|\G|}$, the solution $\bm{\alpha}$ of the optimization task \eqref{eq:J(alfa,z)} is
equal to the solution of \eqref{eq:l}:
\begin{equation}
  \argmin_{\bm{\alpha}\in \CA,
  \bm{z}\in\R_+^{|\G|}} J(\bm{\alpha},\b{z}),\text{ where}\label{eq:min:J(alfa,z)}
\end{equation}
\vspace*{-1cm}
\begin{eqnarray}
      \lefteqn{J(\bm{\alpha},\b{z}) =}\label{eq:J(alfa,z)}\\
      &&\hspace*{-.7cm}{=\frac{1}{2}\left\|\b{x}_{O_t}-(\b{D}_{t-1})_{O_t}\bm{\alpha}\right\|^2_2+\kappa\frac{1}{2}\left(\bm{\alpha}^Tdiag(\bm{\zeta})\bm{\alpha}+\left\|\b{z}\right\|_{\beta}\right)},\nonumber
\end{eqnarray}
$\bm{\zeta}=\bm{\zeta}(\b{z})\in \R^{\da}$ and $\zeta_j=\sum_{G\in\G, G\ni j}\left(d^G_j\right)^2/z^G$.
The optimization of \eqref{eq:J(alfa,z)} can be carried out by iterative alternating steps. One can minimize the quadratic cost function
 on the convex set $\bm{\CA}$ for a given $\b{z}$ with standard solvers \cite{bertsekas99nonlinear}. Then, one can use the
variation principle and find solution $\b{z}=(z^G)_{G\in\G}$ for a fixed $\bm{\alpha}$ by means of the explicit expression
\begin{equation}
    z^G=\|\b{d}^G\circ\bm{\alpha}\|_2^{2-\eta}\|(\|\b{d}^G\circ\bm{\alpha}\|_2)_{G\in\G}\|^{\eta-1}_{\eta}.
\end{equation}
Note that for numerical stability, smoothing $\b{z}=\max(\b{z},\varepsilon)$ ($0<\varepsilon\ll 1$) is suggested in practice.

\subsubsection{Dictionary optimization ($\b{D}$)}
The block-coordinate descent (BCD) method \cite{bertsekas99nonlinear} is used for the optimization of $\b{D}$: columns $\b{d}_j$ in $\b{D}$ are optimized
one-by-one by keeping the other columns
($\b{d}_i,i\ne j$) fixed. For a given $j$, $\hat{f}_t$ is quadratic in $\b{d}_j$. The minimum is found by
solving $\frac{\partial \hat{f}_t}{\partial \b{d}_j}(\b{u}_j)=\b{0}$, and then this solution is projected to the constraint
set $\CD_j$ ($\b{d}_j\leftarrow \Pi_{\CD_j}(\b{u}_j)$). One can show by executing the differentiation that $\b{u}_j$ satisfies the linear equation system
\begin{equation}
  \b{C}_{j,t}\b{u}_j=\b{b}_{j,t}-\b{e}_{j,t}+\b{C}_{j,t}\b{d}_j,
\end{equation}
where
\begin{align}
\b{C}_{j,t} &= \sum_{i=1}^t\left(\frac{i}{t}\right)^{\rho}\D_i \alpha_{i,j}^2\in\R^{\dx\times \dx},\\
\b{e}_{j,t} &= \sum_{i=1}^t\left(\frac{i}{t}\right)^{\rho}\D_i\b{D}\bm{\alpha}_{i}\alpha_{i,j}\in\R^{\dx},\\
\b{B}_t &= \sum_{i=1}^t\left(\frac{i}{t}\right)^{\rho}\D_i\b{x}_i\bm{\alpha}_{i}^T=[\b{b}_{1,t},\ldots,\b{b}_{\da,t}],
\end{align}
matrices $\b{C}_{j,t}$ are diagonal, $\b{B}_t\in\R^{\dx\times\da}$, and $\D_i\in\R^{\dx\times\dx}$ is the diagonal matrix representation of the $O_i$ set (for $j\in O_i$  the $j^{th}$ diagonal is 1 and is $0$ otherwise). It is sufficient to
update statistics $\{\{\b{C}_{j,t}\}_{j=1}^{\da},\b{B}_t,\{\b{e}_{j,t}\}_{j=1}^{\da}\}$ online for the optimization of $\hat{f}_t$, which can be done exactly for $\b{C}_{j,t}$ and $\b{B}_{t}$:
\begin{align}
 \b{C}_{j,t} & =\gamma_t\b{C}_{j,t-1}+\D_{t}\alpha_{tj}^2,\\
 \b{B}_{t} &= \gamma_t\b{B}_{t-1}+\D_{t}\b{x}_{t}\bm{\alpha}_{t}^T,
\end{align}
where $\gamma_t=\left(1-\frac{1}{t}\right)^{\rho}$ and the recursions are initialized by (i) $\b{C}_{j,0}=\b{0}$, $\b{B}_{0}=\b{0}$ for $\rho=0$ and (ii) in an arbitrary way for $\rho>0$. According to
numerical experiences,
\begin{equation}
\b{e}_{j,t} =  \gamma_t \b{e}_{j,t-1}+\D_t\b{D}_t\bm{\alpha}_{t}\alpha_{t,j}, \label{eq:e:general}
\end{equation}
is a good approximation for $\b{e}_{j,t}$ with the actual estimation $\b{D}_t$ and with initialization $\b{e}_{j,0}=\b{0}$. It may be worth noting that the
convergence speed is often improved if statistics are updated in mini-batches  $\{\b{x}_{O_{t,1}},\ldots,\b{x}_{O_{t,R}}\}$.\footnote{The Matlab code of the OSDL method is available at \url{http://nipg.inf.elte.hu/szzoli}.} 

\section{OSDL Based Collaborative Filtering}\label{sec:OSDL via CF}
We formulate the CF task as an OSDL optimization problem in Section~\ref{sec:CF casted as OSDL}. According to the CF literature,
oftentimes neighbor-based corrections improve the precision of the estimation. We also use this technique (Section~\ref{sec:NN correction}) to improve the
OSDL estimations.

\subsection{CF Casted as an OSDL Problem}\label{sec:CF casted as OSDL}
Below, we transform the CF task into an OSDL problem. Consider the  $t^{th}$ user's known ratings as OSDL observations $\b{x}_{O_t}$. Let
the optimized group-structured dictionary on these observations be $\b{D}$. Now, assume that we have a test user and
his/her ratings, i.e., $\b{x}_O\in\R^{|O|}$. The task is to estimate $\b{x}_{\{1,\ldots,\dx\}\backslash O}$, that is, the
missing coordinates of $\b{x}$ (the missing ratings of the user) that can be accomplished as follows:
\begin{enumerate}
  \item
      Remove the rows of the non-observed ${\{1,\ldots,\dx\}\backslash O}$ coordinates from $\b{D}$. The obtained $|O|\times \da$ sized matrix $\b{D}_O$ and $\b{x}_O$ can be used to estimate $\bm{\alpha}$ by solving
      the structured sparse coding problem \eqref{eq:l}.
  \item Using the estimated representation $\bm{\alpha}$, estimate $\b{x}$ as
      \begin{equation}
    \hat{\b{x}}=\b{D}\bm{\alpha}.
      \end{equation}
\end{enumerate}

\subsection{Neighbor Based Correction}\label{sec:NN correction}
According to the CF literature, neighbor based correction schemes may further improve the precision of the estimations \cite{ricci11recommender}.
This neighbor correction approach
\begin{itemize}
  \item relies on the assumption that similar items (e.g., jokes/movies) are rated similarly and
  \item can be adapted to OSDL-based CF estimation in a natural fashion.
\end{itemize}
Here, we detail the idea. Let us assume that the similarities $s_{ij}\in\R$ ($i,j\in\{1,\ldots,\dx\}$) between individual items are given.
We shall provide similarity forms in Section~\ref{sec:item similarities}. Let  $\b{d}_k\bm{\alpha}_t\in\R$ be the OSDL estimation
for the rating of the $k^{th}$ non-observed item of the $t^{th}$ user ($k\not\in O_t$), where $\b{d}_k\in\R^{1\times \da}$ is
the $k^{th}$ row of matrix $\b{D}\in\R^{\dx\times\da}$, and $\bm{\alpha}_t\in\R^{\da}$ is computed according to
Section~\ref{sec:CF casted as OSDL}.

Let the prediction error on the observable item neighbors ($j$) of the $k^{th}$ item of the $t^{th}$ user ($j\in O_t\backslash \{k\}$) be
$\b{d}_j\bm{\alpha}_t-x_{jt}\in\R$. These prediction errors can be used for the correction of the OSDL estimation ($\b{d}_k\bm{\alpha}_t$) by taking into account the  $s_{ij}$ similarities:
\begin{align}
    \hat{x}_{kt} &= \b{d}_k\bm{\alpha}_t + \gamma_1 \left[\frac{\sum_{j\in O_t\backslash \{k\}}s_{kj}(\b{d}_j\bm{\alpha}_t-x_{jt})}{\sum_{j\in O_t\backslash \{k\}}s_{kj}}\right],\text{ or}\label{eq:xhat}\\
    \hat{x}_{kt} &= \gamma_0(\b{d}_k\bm{\alpha}_t) + \gamma_1 \left[\frac{\sum_{j\in O_t\backslash \{k\}}s_{kj}(\b{d}_j\bm{\alpha}_t-x_{jt})}{\sum_{j\in O_t\backslash \{k\}}s_{kj}}\right],\label{eq:xhat-with-0}
\end{align}
where $k\not\in O_t$. Here, \eqref{eq:xhat} is analogous to the form of \cite{takacs09scalable}, \eqref{eq:xhat-with-0} is a simple modification: it modulates the first term with a separate $\gamma_0$ weight.

\section{Numerical Results}\label{sec:numerical results}
We have chosen the Jester dataset (Section~\ref{sec:Jester}) for the illustration of the OSDL based CF approach.
It is a standard benchmark for CF. We detail our preferred item similarities in Section~\ref{sec:item similarities}.
To evaluate the CF based estimation, we use the performance measures given in Section~\ref{sec:performance measure}.
Section~\ref{sec:evaluation} is about our numerical experiences.

\subsection{The Jester Dataset}\label{sec:Jester}
The dataset \cite{goldber01eigentaste} contains $4,136,360$ ratings from $73,421$ users to $100$ jokes on a
continuous $[-10,10]$ range. The worst and best possible gradings are $-10$ and $+10$, respectively. A fixed $10$ element subset of
the jokes is called gauge set and it was evaluated by all users. Two third of the users have rated at least $36$ jokes, and the remaining ones have rated between $15$ and $35$ jokes. The average number of user ratings per joke is $46$.

\subsection{Item Similarities}\label{sec:item similarities}
In the neighbor correction step \eqref{eq:xhat} or \eqref{eq:xhat-with-0} we need the $s_{ij}$ values representing the
similarities of the $i^{th}$ and $j^{th}$ items. We define this value as the similarity of the $i^{th}$ and $j^{th}$
rows ($\b{d}_i$ and $\b{d}_j$) of the optimized OSDL dictionary $\b{D}$ \cite{takacs09scalable}:
    \begin{align}
      S_1 :\quad & s_{ij}=s_{ij}(\b{d}_i,\b{d}_j)  = \left(\frac{\max(0,\left<\b{d}_i,\b{d}_j\right>)}{\left\|\b{d}_i\right\|_2\left\|\b{d}_j\right\|_2}\right)^{\beta}\text{, or}\label{eq:S1}\\
      S_2 :\quad & s_{ij} = s_{ij}(\b{d}_i,\b{d}_j) = \left(\frac{\left\|\b{d}_i-\b{d}_j\right\|_2^2}{\left\|\b{d}_i\right\|_2\left\|\b{d}_j\right\|_2}\right)^{-\beta},\label{eq:S2}
  \end{align}
where $\beta>0$ is the parameter of the similarity measure. Quantities $s_{ij}$ are non-negative; if the value of $s_{ij}$ is close to zero (large) then the $i^{th}$ and $j^{th}$ items are very different (very similar).

\subsection{Performance Measure}\label{sec:performance measure}
In our numerical experiments we used the RMSE (root mean square error) and the MAE (mean absolute error) measure
for the evaluation of the quality of the estimation, since these are the most popular measures in the CF literature.
The RMSE and MAE measure is the average squared/absolute difference of the true and the estimated rating values, respectively:
\begin{align}
 RMSE &= \sqrt{\frac{1}{|\mathscr{S}|}\sum_{(i,t)\in \mathscr{S}} (x_{it}-\hat{x}_{it})^2},\\
 MAE &= \frac{1}{|\mathscr{S}|}\sum_{(i,t)\in \mathscr{S}} |x_{it}-\hat{x}_{it}|,
\end{align}
where $\mathscr{S}$ denotes either the validation or the test set.

\subsection{Evaluation}\label{sec:evaluation}
Here we illustrate the efficiency of the OSDL-based CF estimation on the Jester dataset (Section~\ref{sec:Jester}) using the
RMSE and MAE performance measures (Section~\ref{sec:performance measure}). We start our discussion with the RMSE results. The MAE performance measure led to similar results; for the sake of completeness we report these results at the end of this section.
To the best of our knowledge, the top results on this database are RMSE = $4.1123$ \cite{takacs08matrix} and RMSE = $4.1229$
\cite{takacs09scalable}. Both works are from the same authors. The method in the first paper is called item neighbor and it
makes use of only neighbor information.
In \cite{takacs09scalable}, the authors used a bridge regression based unstructured dictionary learning model---with a neighbor correction scheme---, they optimized the dictionary
by gradient descent and set $\da$ to 100. These are our performance baselines.

To study the capability of the OSDL approach in CF, we focused on the following issues:
\begin{itemize}
 \item Is structured dictionary $\b{D}$ beneficial for prediction purposes, and how does it compare to the dictionary of classical (unstructured) sparse dictionary?
 \item How does the OSDL parameters and the similarity/neighbor correction applied affect the efficiency of the prediction?
 \item How do different group structures $\G$ fit to the CF task?
\end{itemize}

In our numerical studies we chose the Euclidean unit sphere for $\CD_i=S_2^{\dx}$ ($\forall i$), and $\CA=\R^{\da}$, and no additional weighting was applied ($\b{d}^G=\chi_G$, $\forall G\in\G$, where $\chi$ is the indicator function). We set $\eta$ of the group-structured regularizer $\Omega$ to $0.5$. Group structure $\G$ of vector $\bm{\alpha}$ was realized on
\begin{itemize}
 \item a $d\times d$ toroid ($\da=d^2$) with $|\G|=\da$ applying $r\ge 0$ neighbors to define $\G$. For $r=0$ ($\G=\{\{1\},\ldots,\{\da\}\}$) the classical sparse representation based dictionary is recovered.
 \item a hierarchy with a complete binary tree structure. In this case:
      \begin{itemize}
      \item $|\G|=\da$, and group $G$ of $\alpha_i$ contains the $i^{th}$ node and its descendants on the tree, and
      \item the size of the tree is determined by the number of levels $l$. The dimension of the hidden representation is then $\da=2^l-1$.
      \end{itemize}
\end{itemize}
The size $R$ of mini-batches was set either to $8$, or to $16$ and the forgetting factor $\rho$ was chosen from set
$\{0,\frac{1}{64},\frac{1}{32},\frac{1}{16},\frac{1}{8},\frac{1}{4},\frac{1}{2},1\}$. The $\kappa$ weight of structure inducing regularizer  $\Omega$ was chosen from the set
$\{\frac{1}{2^{-1}},\frac{1}{2^{0}},\frac{1}{2^{1}},\frac{1}{2^{2}},\frac{1}{2^{4}},\frac{1}{2^{6}},\ldots,\frac{1}{2^{14}}\}$.
We studied similarities $S_1$, $S_2$ [see \eqref{eq:S1}-\eqref{eq:S2}] with both neighbor correction schemes
[\eqref{eq:xhat}-\eqref{eq:xhat-with-0}]. In what follows, corrections based on \eqref{eq:xhat} and \eqref{eq:xhat-with-0} will
be called $S_1$, $S_2$ and $S_1^0$, $S_2^0$, respectively. Similarity parameter $\beta$ was chosen from the set
$\{0.2,1,1.8,2.6,\ldots,14.6\}$. In the BCD step of the optimization of $\b{D}$, $5$ iterations were applied. In the
$\bm{\alpha}$ optimization step, we used $5$ iterations, whereas smoothing parameter $\epsilon$ was $10^{-5}$.

We used a $90\%-10\%$ random split for the observable ratings in our experiments, similarly to \cite{takacs09scalable}:
\begin{itemize}
 \item training set ($90\%$) was further divided into 2 parts:
    \begin{itemize}
    \item we chose the $80\%$ observation set $\{O_t\}$ randomly, and optimized $\b{D}$ according to the corresponding $\b{x}_{O_t}$ observations,
    \item we used the remaining $10\%$ for validation, that is for choosing the optimal
	  OSDL parameters ($r$ or $l$, $\kappa$, $\rho$), BCD optimization parameter ($R$),
	  neighbor correction ($S_1$, $S_2$, $S_1^0$, $S_2^0$), similarity parameter ($\beta$), and
	  correction weights ($\gamma_i$s in \eqref{eq:xhat} or \eqref{eq:xhat-with-0}).
    \end{itemize}
 \item we used the remaining $10\%$ of the data for testing.
\end{itemize}
The optimal parameters were estimated on the validation set, and then used on the test set. The resulting RMSE/MAE score was the performance of the estimation.

\subsubsection{Toroid Group Structure.}
In this section we provide results using toroid group structure. We set $d=10$. The size of the toroid was $10\times 10$, and thus the dimension of the representation was $\da=100$.

\begin{figure*}%
\centering%
\subfloat[][]{\includegraphics[width=6.05cm]{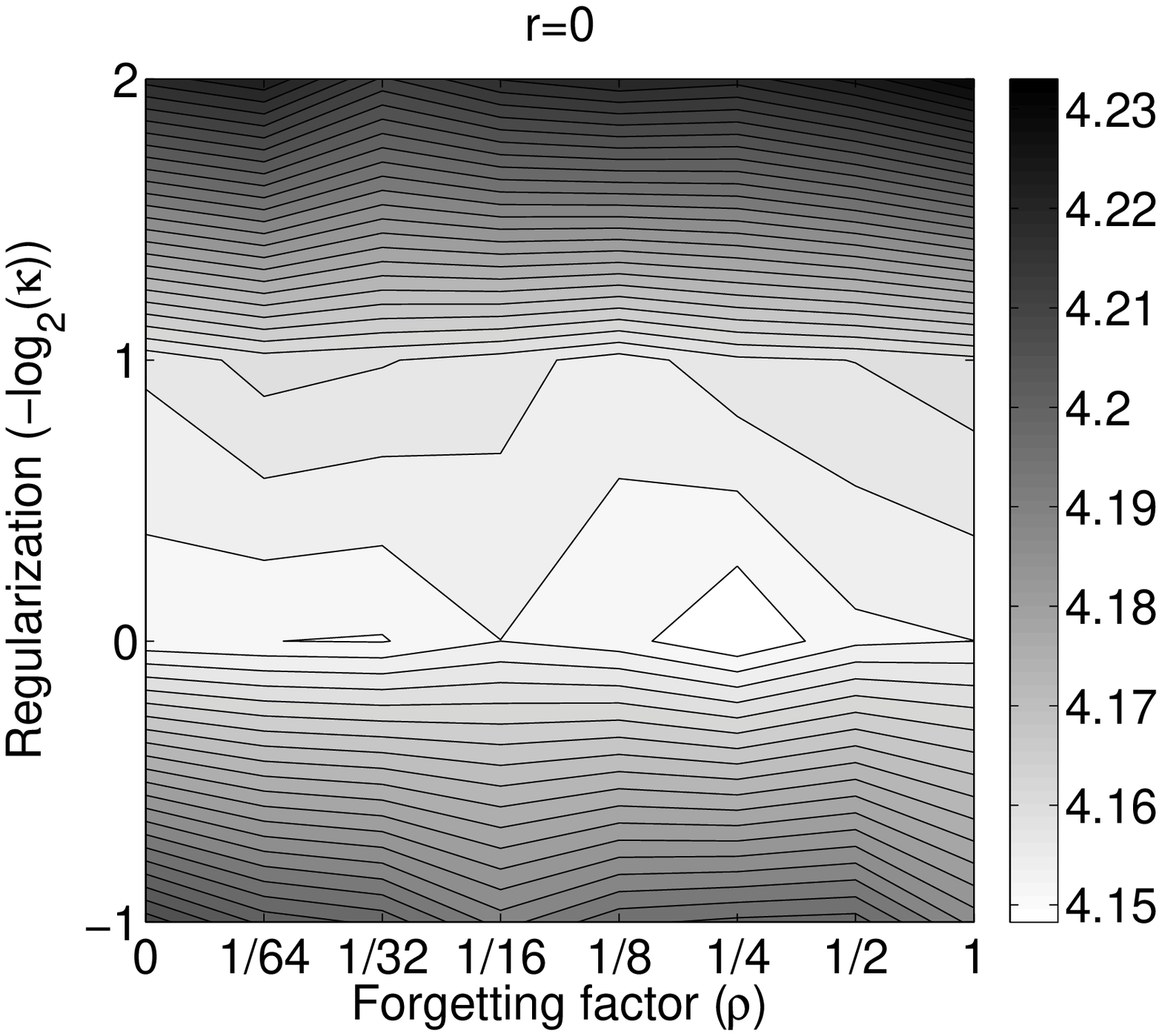}}
\subfloat[][]{\includegraphics[width=6.05cm]{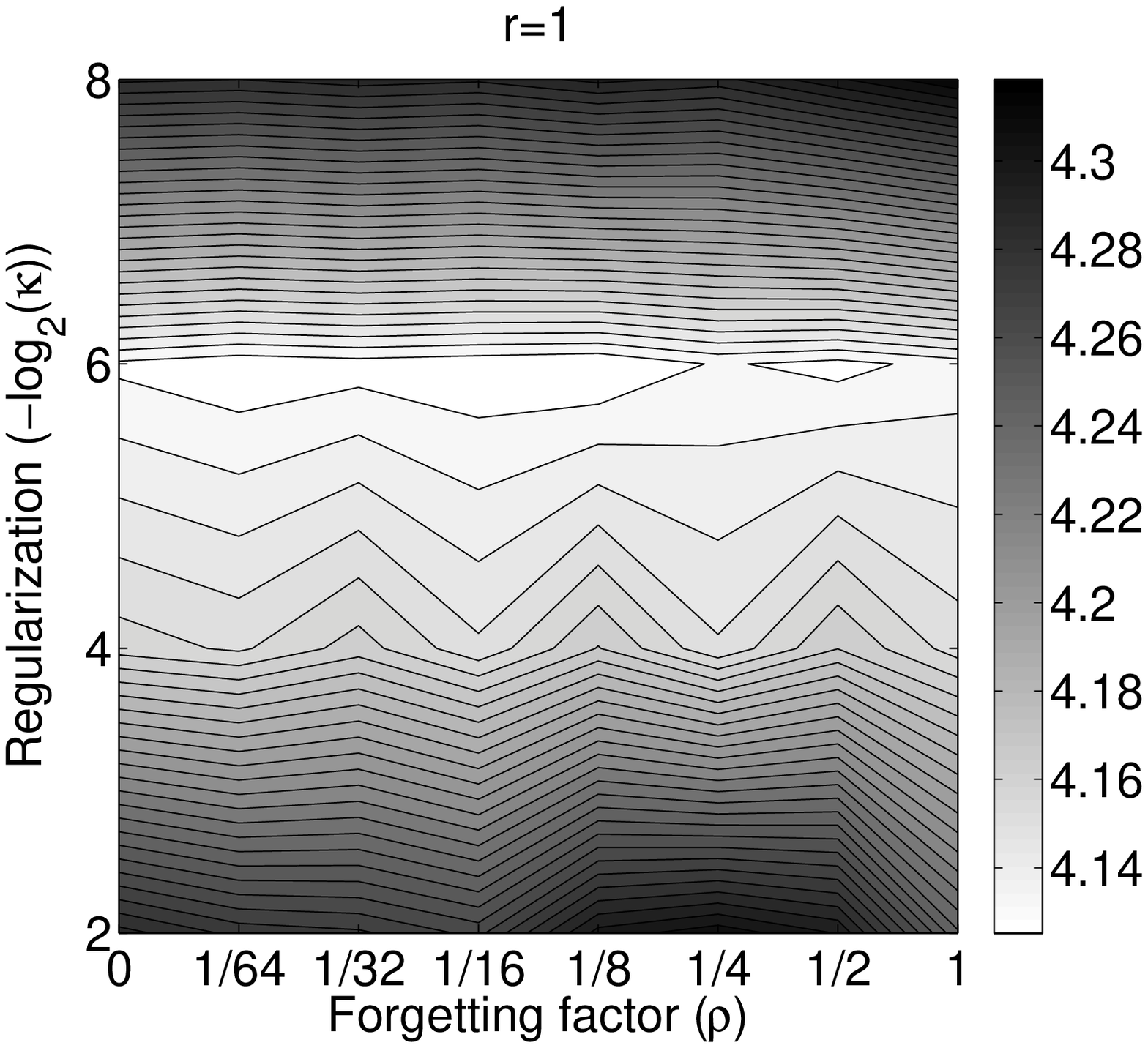}}
\subfloat[][]{\includegraphics[width=6.05cm]{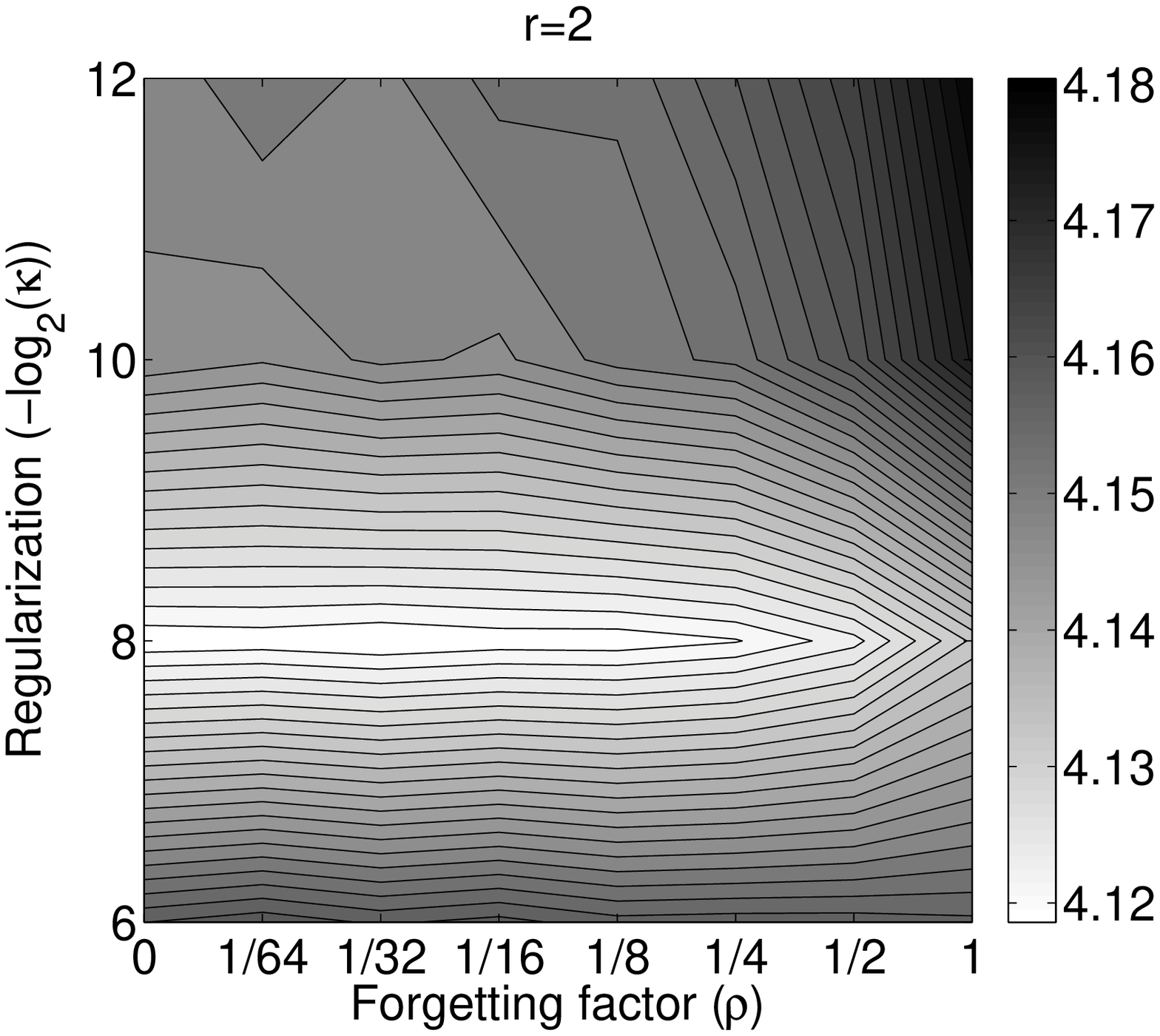}}\\
\subfloat[][]{\includegraphics[width=6.05cm]{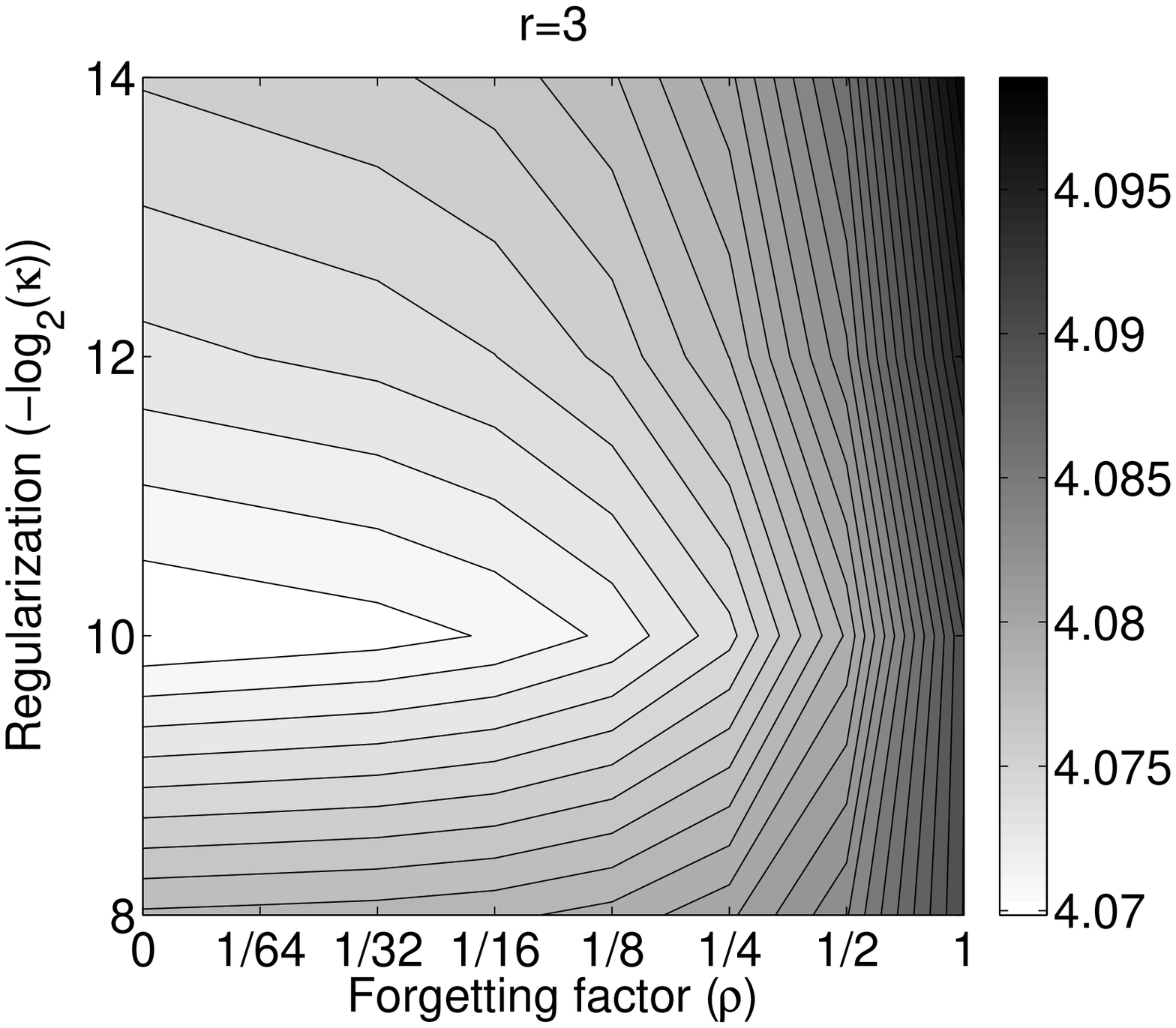}}
\subfloat[][]{\includegraphics[width=6.05cm]{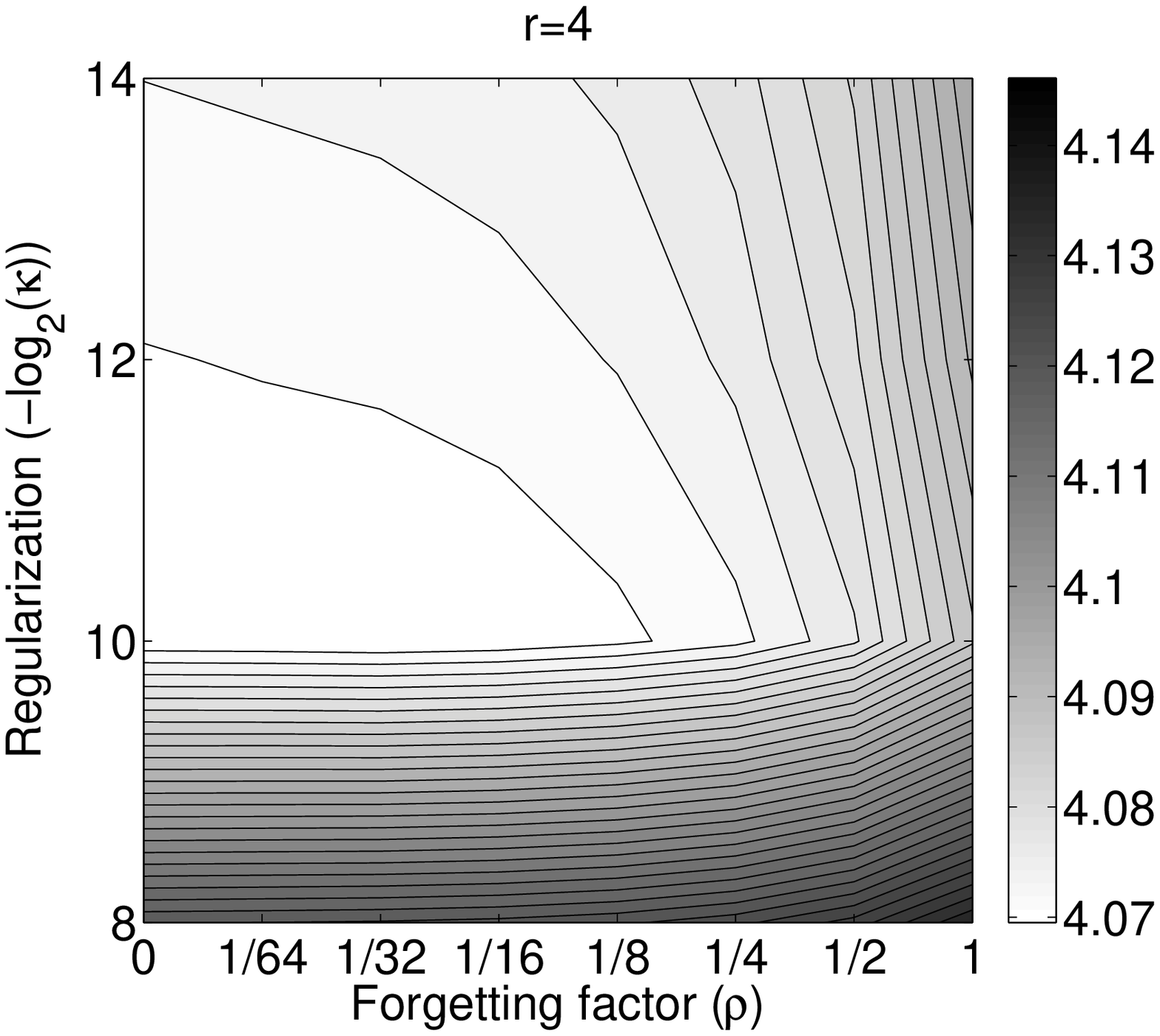}}
\subfloat[][]{\includegraphics[width=6.05cm]{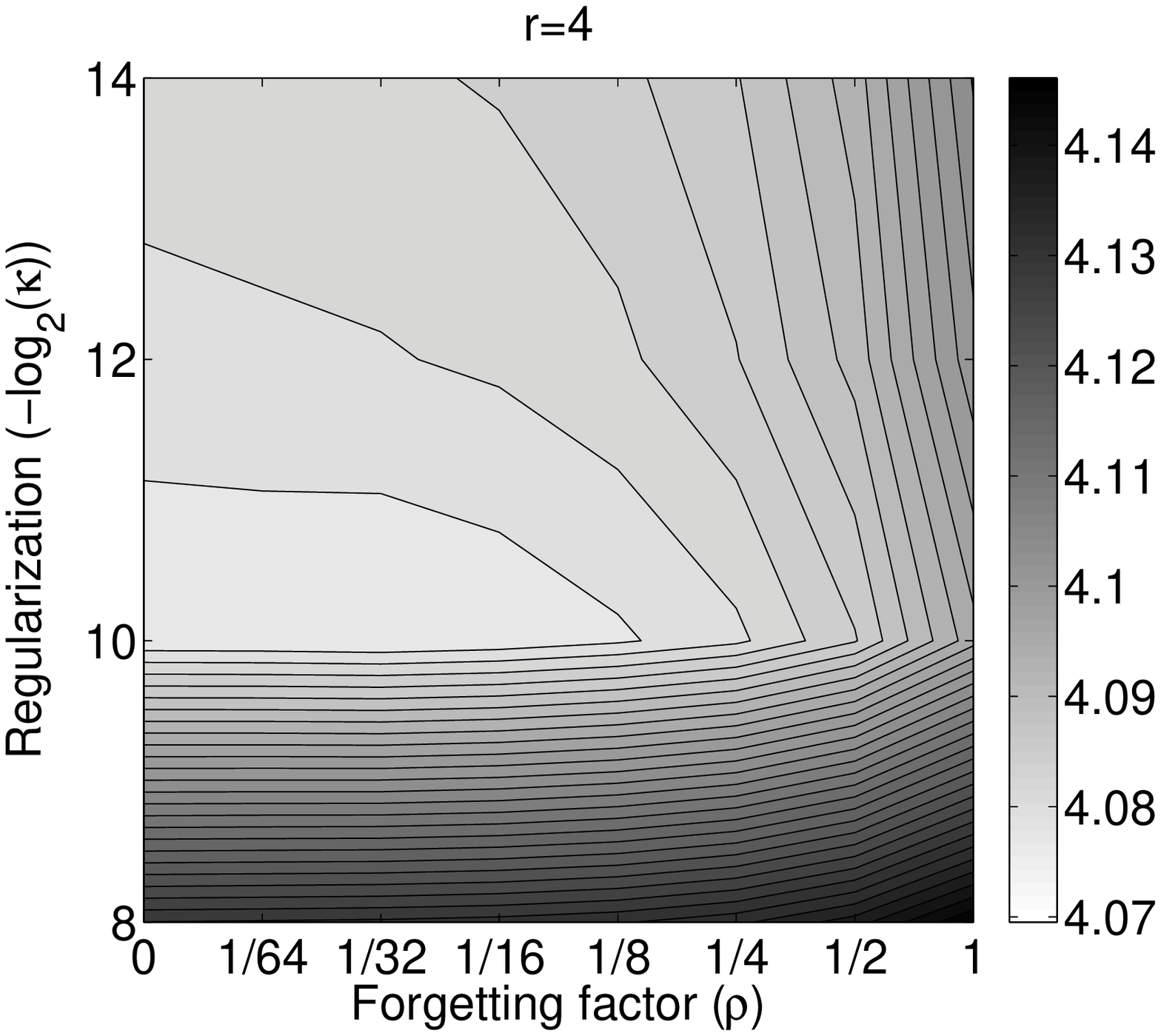}}
\caption[]{RMSE validation surfaces [(a)-(e)] and test surfaces  (f) as a function of forgetting factor ($\rho$) and regularization ($\kappa$).
 For a fixed $(\kappa, \rho)$ parameter pair, the surfaces show the best RMSE values optimized in the $\beta$ similarity parameter. The group structure ($\G$) is  toroid. The applied neighbor correction was $S_1^0$. (a): $r=0$ (no structure). (b): $r=1$. (c): $r=2$. (d): $r=3$. (e)-(f): $r=4$, on the same scale.}%
\label{fig:torus:validation surfaces}
\end{figure*}

\begin{figure}[h]%
\centering%
\includegraphics[width=6.2cm]{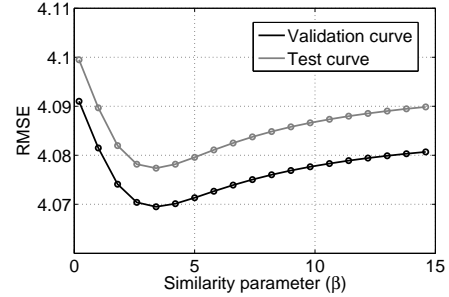}
\caption[]{RMSE validation and test curves for toroid group structure using the optimal neighborhood size $r=4$, regularization weight $\kappa=\frac{1}{2^{10}}$, forgetting factor $\rho=\frac{1}{2^5}$, mini-batch size $R=8$, and similarity parameter $\beta=3.4$. The applied neighbor correction was $S_1^0$.}%
\label{fig:torus:validation vs test curve}
\end{figure}

In the \textbf{first experiment} we study how the size of neighborhood ($r$) affects the results.
This parameter corresponds to the ``smoothness'' imposed on the
group structure: when $r=0$, then there is no relation between the $\b{d}^j\in\R^{\da}$ columns in $\b{D}$ 
(no structure). As we increase  $r$, the $\b{d}^j$ feature vectors will be more and more aligned in a smooth way. To this end, we set the neighborhood size to $r=0$ (no structure), and then increased it to $1$, $2$, $3$, $4$, and $5$. For each $(\kappa,\rho,\beta)$, we calculated the RMSE of our
estimation, and then for each fixed  ($\kappa,\rho$) pair, we minimized these RMSE values in $\beta$. The resulting validation and test surfaces are shown in Fig.~\ref{fig:torus:validation surfaces}.
For the best ($\kappa,\rho$) pair, we also present the RMSE values as a function of $\beta$  (Fig.~\ref{fig:torus:validation vs test curve}).
In this illustration we used $S_1^0$ neighbor correction and $R=8$ mini-batch size. We note that we got similar results using $R=16$ too. Our results can be summarized as follows.

\begin{itemize}
  \item For a fixed neighborhood parameter $r$, we have that:
      \begin{itemize}
    \item The validation and test surfaces are very similar (see Fig.~\ref{fig:torus:validation surfaces}(e)-(f)).
          It implies that the validation surfaces are good indicators for the test errors. For the best $r$, $\kappa$ and $\rho$ parameters, we can observe that the validation and test curves (as functions of $\beta$) are very similar. This is demonstrated in Fig.~\ref{fig:torus:validation vs test curve}, where we used $r=4$ neighborhood size and $S_1^0$ neighbor correction. We can also notice that
          (i) both curves have only one local minimum, and (ii) these minimum points are close to each other.
    \item The quality of the estimation depends mostly on the $\kappa$ regularization parameter. As we increase $r$, the best $\kappa$ value is decreasing.
    \item The estimation is robust to the different choices of forgetting factors (see Fig.~\ref{fig:torus:validation surfaces}(a)-(e)). In other words, this parameter $\rho$ can help in fine-tuning the results.
      \end{itemize}
  \item Structured dictionaries ($r>0$) are advantageous over those methods that do not impose structure on the dictionary elements ($r=0$). For $S_1^0$ and $S_2^0$ neighbor corrections, we summarize the RMSE results in
    Table~\ref{tab:torus:perf:r}. Based on this table we can conclude that in the studied parameter domain
         \begin{itemize}
        \item the estimation is robust to the selection of the mini-batch size ($R$). We got the best results using $R=8$. Similarly to the role of parameter $\rho$, adjusting $R$ can be used for fine-tuning.
        \item the $S_1^0$ neighbor correction lead to the smallest RMSE value.
        \item When we increase $r$ up to $r=4$,  the results improve. However, for $r=5$, the RMSE values do not improve anymore; they are about the same that we have using $r=4$.
        \item The smallest RMSE we could achieve was $4.0774$, and the best known result so far was
              RMSE = $4.1123$ \cite{takacs08matrix}. This proves the efficiency of our OSDL based collaborative filtering algorithm.
        \item We note that our RMSE result seems to be significantly better than the that of the competitors: we repeated this experiment $5$ more times with different randomly selected training, test, and validation sets, and our RMSE results have never been worse than $4.08$.
          \end{itemize}
\end{itemize}

In the \textbf{second experiment} we studied how the different neighbor corrections ($S_1$, $S_2$, $S_1^0$, $S_2^0$) affect the performance of the proposed algorithm. To this end, we set the neighborhood parameter to $r=4$ because it proved to be optimal in the previous experiment. Our results are summarized in Table~\ref{tab:torus:perf:S}. From these results we can observe that
\begin{itemize}
    \item our method is robust to the selection of correction methods. Similarly to the $\rho$ and $R$ parameters, the
    neighbor correction scheme can help in fine-tuning the results.
    \item The introduction of $\gamma_0$ in \eqref{eq:xhat-with-0} with the application of $S_1^0$ and $S_2^0$ instead of
    $S_1$ and $S_2$ proved to be advantageous in the neighbor correction phase.
    \item For the studied CF problem, the $S_1^0$ neighbor correction method (with $R=8$) lead to the smallest RMSE value, $4.0774$.
    \item The $R\in\{8,16\}$ setting yielded us similarly good results. Even with $R=16$, the RMSE value was $4.0777$.
\end{itemize}

\begin{table}
    \centering
    \caption{Performance (RMSE) of the OSDL prediction using toroid group structure ($\G$) with different neighbor sizes $r$ ($r=0$: unstructured case). First-second row: mini-batch size $R=8$, third-fourth row: $R=16$.
        Odd rows: $S_1^0$, even rows: $S_2^0$ neighbor correction.  For fixed $R$, the best performance is highlighted with boldface typesetting.}
    \begin{tabular}{|c|c|c|c|c|c|c|}
    \hline
       &       &  $r=0$ & $r=1$ & $r=2$ & $r=3$ & $r=4$ \\
    \hline\hline
     $R=8$ & $S_1^0$ & $4.1594$ & $4.1326$ & $4.1274$ & $4.0792$ &  $\mathbf{4.0774}$\\
     \hline
           & $S_2^0$ & $4.1765$ & $4.1496$ & $4.1374$ &  $4.0815$&  $4.0802$\\
     \hline\hline
     $R=16$ & $S_1^0$ & $4.1611$ & $4.1321$ & $4.1255$ & $4.0804$ & $\mathbf{4.0777}$\\
    \hline
            & $S_2^0$ & $4.1797$ & $4.1487$ & $4.1367$ & $4.0826$ &  $4.0802$\\
    \hline
    \end{tabular}
    \label{tab:torus:perf:r}
\end{table}

\begin{table}
    \caption{Performance (RMSE) of the OSDL prediction for different neighbor corrections
    using toroid group structure ($\G$). Columns: applied neighbor corrections. Rows: mini-batch size $R=8$ and $16$.
     The neighbor size was set to $r=4$. For fixed $R$, the best performance is highlighted with boldface typesetting.}
    \centering
    \begin{tabular}{|c|c|c|c|c|}
    \hline
          & $S_1$ & $S_2$ & $S_1^0$ & $S_2^0$\\
   \hline\hline
    $R=8$ & $4.0805$ & $4.0844$ & $\mathbf{4.0774}$ & $4.0802$\\
   \hline
    $R=16$ & $4.0809$ & $4.0843$ & $\mathbf{4.0777}$ & $4.0802$\\
    \hline
    \end{tabular}
    \label{tab:torus:perf:S}
\end{table}

\subsubsection{Hierarchical Group Structure.}
In this section we provide results using hierarchical $\bm{\alpha}$ representation. The group structure $\G$ was chosen to represent a complete binary tree.

In our \textbf{third experiment} we study how the number of levels ($l$) of the tree affects the results. To this end, we set the number of levels to $l=3$, $4$, $5$, and $6$.
Since $\da$, the dimension of the hidden representation $\bm{\alpha}$, equals to $2^l-1$, these
$l$ values give rise to dimensions $\da=7$, $15$, $31$, and $63$. Validation and test surfaces are provided in Fig.~\ref{fig:bintree:validation surfaces,val-vs-test
curve}(a)-(c) and (e)-(f), respectively.
The surfaces show for each $(\kappa,\rho)$ pair, the minimum RMSE values taken in the similarity parameter $\beta$.
For the best $(\kappa,\rho)$ parameter pair, the dependence of RMSE on
$\beta$ is presented in Fig.~\ref{fig:bintree:validation surfaces,val-vs-test curve}(d).
In this illustration we used $S_1^0$ neighbor correction, and the mini-batch size was set to $R=8$. Our results are summarized below. We note that we obtained similar results with mini-batch size $R=16$.

\begin{itemize}
  \item For fixed number of levels $l$, similarly to the toroid group structure (where the size $r$ of the neighborhood was fixed),
      \begin{itemize}
      \item validation and test surfaces are very similar, see
      Fig.~\ref{fig:bintree:validation surfaces,val-vs-test curve}(b)-(c).
      Validation and test curves as a function of $\beta$ behave alike, see Fig.~\ref{fig:bintree:validation surfaces,val-vs-test curve}(d).
      \item the precision of the estimation depends mostly on the regularization parameter
      $\kappa$; forgetting factor $\rho$ enables fine-tuning.
      \end{itemize}
  \item The obtained RMSE values are summarized in Table~\ref{tab:bintree:perf:l} for
      $S_1^0$ and $S_2^0$ neighbor corrections. According to the table, the quality of estimation
      is about the same for mini-batch size $R=8$ and $R=16$; the $R=8$ based estimation seems somewhat more precise.
      Considering the neighbor correction schemes $S_1^0$ and $S_2^0$, $S_1^0$ provided better predictions.
  \item As a function of the number of levels, we got the best result for $l=4$, RMSE = $4.1220$; RMSE values decrease
      until $l=4$ and then increase for $l>4$.
  \item Our best obtained RMSE value is $4.1220$; it was achieved for dimension only $\da=15$. We note that this small dimensional,
      hierarchical group structure based result is also better than that of \cite{takacs09scalable} with RMSE = $4.1229$, which makes use of unstructured dictionaries with $\da=100$.
     The result is also competitive with the RMSE = $4.1123$ value of \cite{takacs08matrix}.
\end{itemize}

In our \textbf{fourth experiment} we investigate how the different neighbor corrections ($S_1$, $S_2$, $S_1^0$,
$S_2^0$) affect the precision of the estimations. We fixed the number of levels to $l=4$, since it proved to be the optimal choice in our
previous experiment. Our results are summarized in Table~\ref{tab:bintree:perf:S}. We found that
\begin{itemize}
  \item the estimation is robust to the choice of neighbor corrections,
  \item it is worth including weight $\gamma_0$ [see \eqref{eq:xhat-with-0}] to improve the precision of prediction,
    that is, to apply correction $S_1^0$ and $S_2^0$ instead of $S_1$ and $S_2$, respectively.
  \item the studied $R\in\{8,16\}$ mini-batch sizes provided similarly good results.
  \item for the studied CF problem the best RMSE value was achieved using $S_1^0$ neighbor correction and mini-batch size $R=8$.
\end{itemize}

\begin{table}
    \centering
    \caption{Performance (RMSE) of the OSDL prediction  for different number of levels ($l$) using binary tree structure ($\G$). First-second row: mini-batch size $R=8$, third-fourth row: $R=16$.
        Odd rows: $S_1^0$, even rows: $S_2^0$ neighbor correction. For fixed $R$, the best performance is highlighted with boldface typesetting.}
    \begin{tabular}{|c|c|c|c|c|c|c|}
    \hline
       &       &  $l=3$ & $l=4$ & $l=5$ & $l=6$ \\
    \hline\hline
     $R=8$ & $S_1^0$ & $4.1572$ & $\mathbf{4.1220}$ & $4.1241$ & $4.1374$\\
     \hline
           & $S_2^0$ & $4.1669$ & $4.1285$ & $4.1298$ &  $4.1362$\\
     \hline\hline
     $R=16$ & $S_1^0$ & $4.1578$ & $4.1261$ & $\mathbf{4.1249}$ & $4.1373$\\
    \hline
            & $S_2^0$ & $4.1638$ & $4.1332$ & $4.1303$ & $4.1383$\\
    \hline
    \end{tabular}
    \label{tab:bintree:perf:l}
\end{table}

\begin{table}
    \caption{Performance (RMSE) of the OSDL prediction for different neighbor corrections using binary tree structure ($\G$).
    Rows: mini-batch size $R=8$ and $16$. Columns: neighbor corrections. Neighbor size: $r=4$. For fixed $R$, the best performance is highlighted with boldface typesetting.}
    \centering
    \begin{tabular}{|c|c|c|c|c|}
    \hline
          & $S_1$ & $S_2$ & $S_1^0$ & $S_2^0$\\
   \hline\hline
    $R=8$ & $4.1255$ & $4.1338$ & $\mathbf{4.1220}$ & $4.1285$\\
   \hline
    $R=16$ & $4.1296$ & $4.1378$ & $\mathbf{4.1261}$ & $4.1332$\\
    \hline
    \end{tabular}
    \label{tab:bintree:perf:S}
\end{table}

\begin{figure*}%
\centering%
\hspace*{0.2cm}\subfloat[][]{\includegraphics[width=6.05cm]{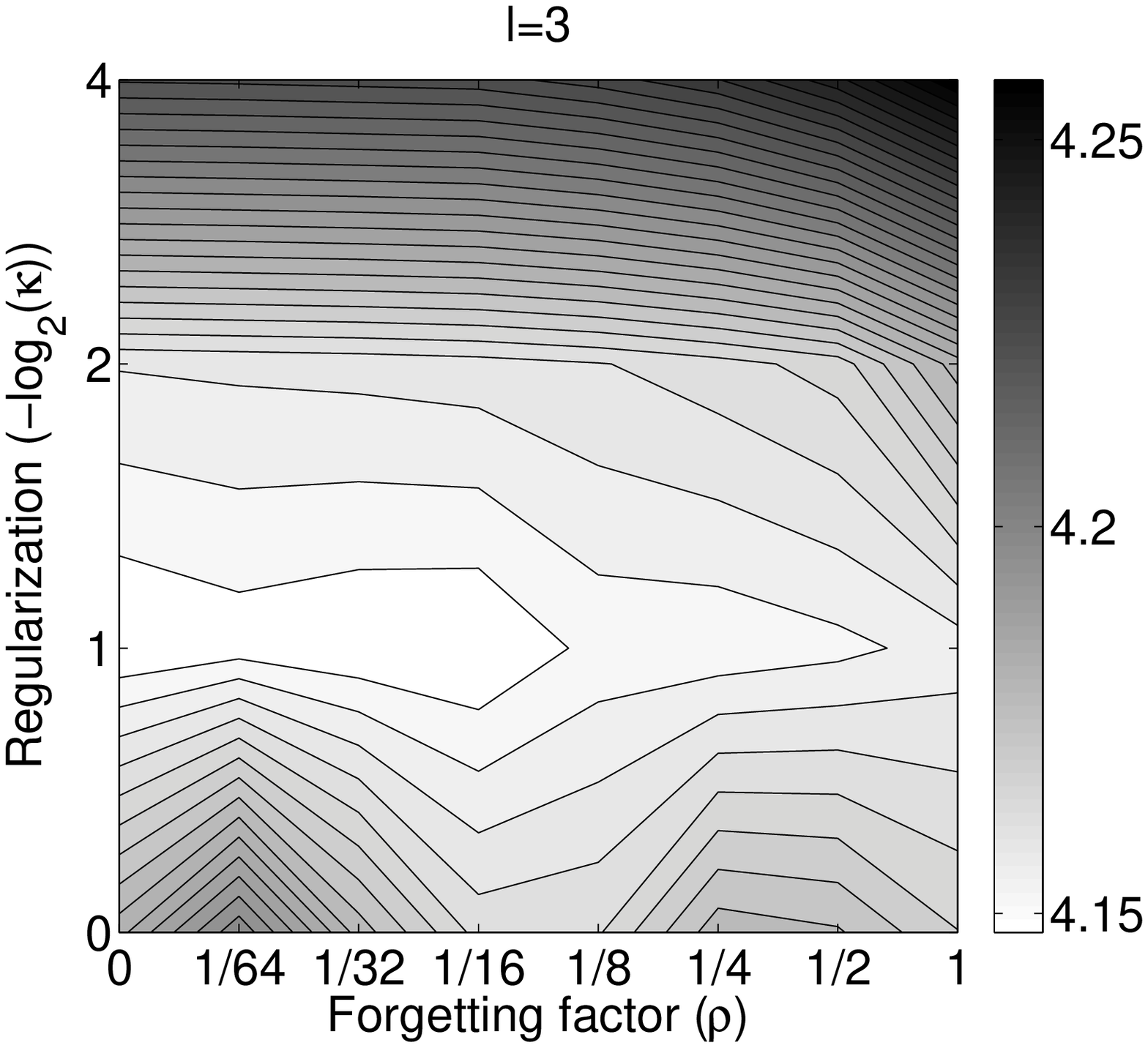}}
\subfloat[][]{\includegraphics[width=6.05cm]{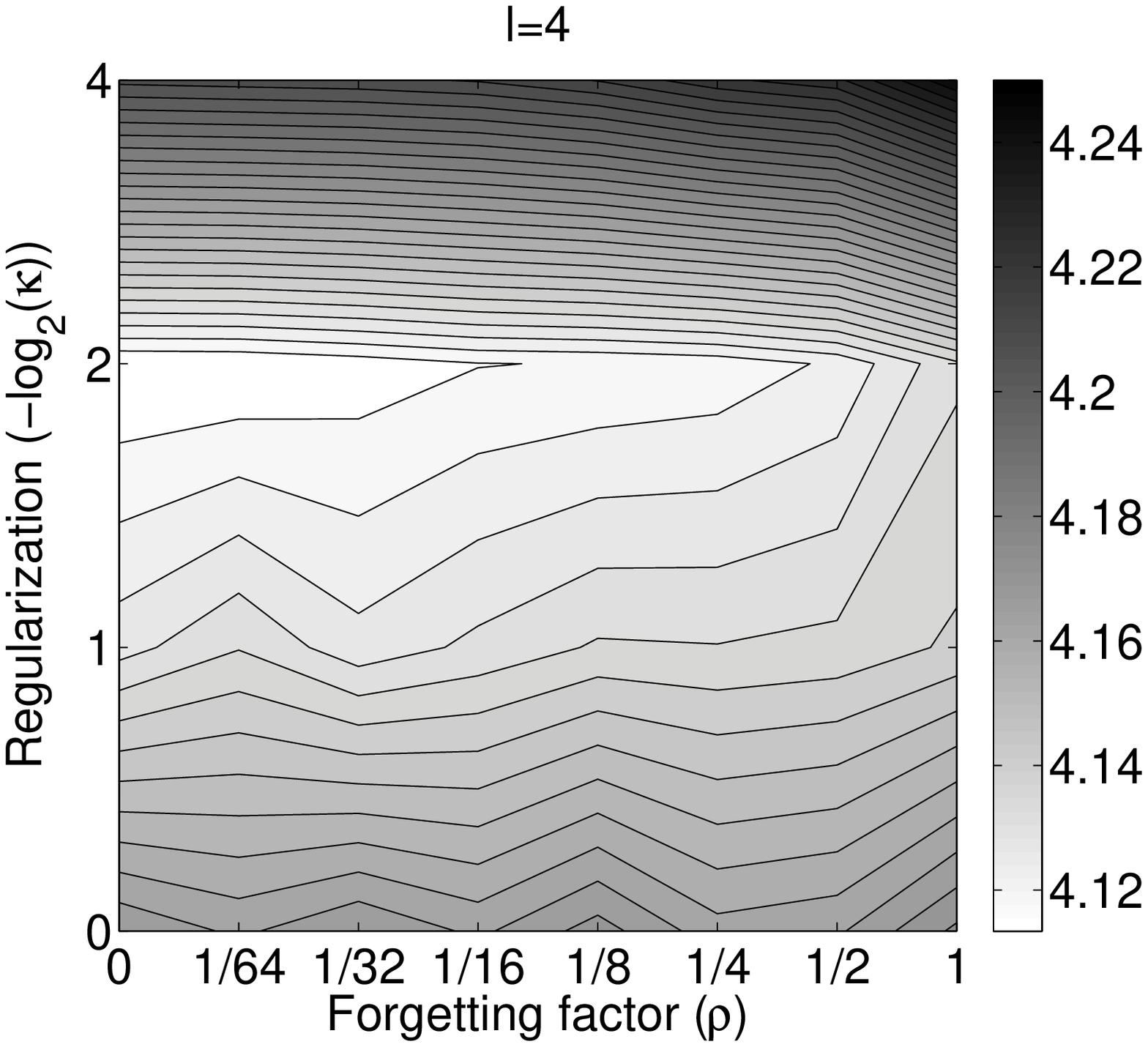}}
\subfloat[][]{\includegraphics[width=6.05cm]{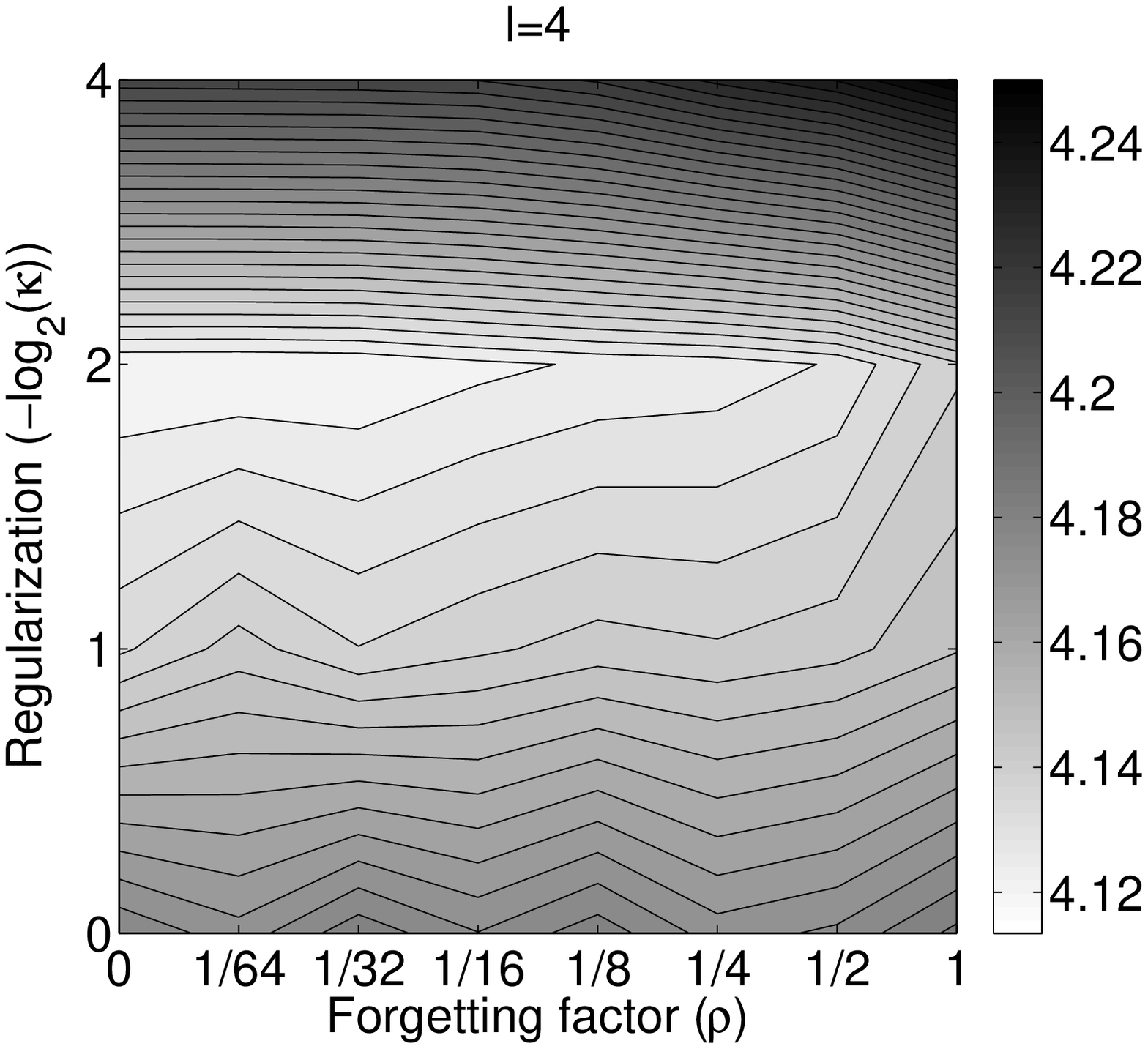}}\\
\subfloat[][]{\includegraphics[width=6.2cm,trim=15 -30 0 0]{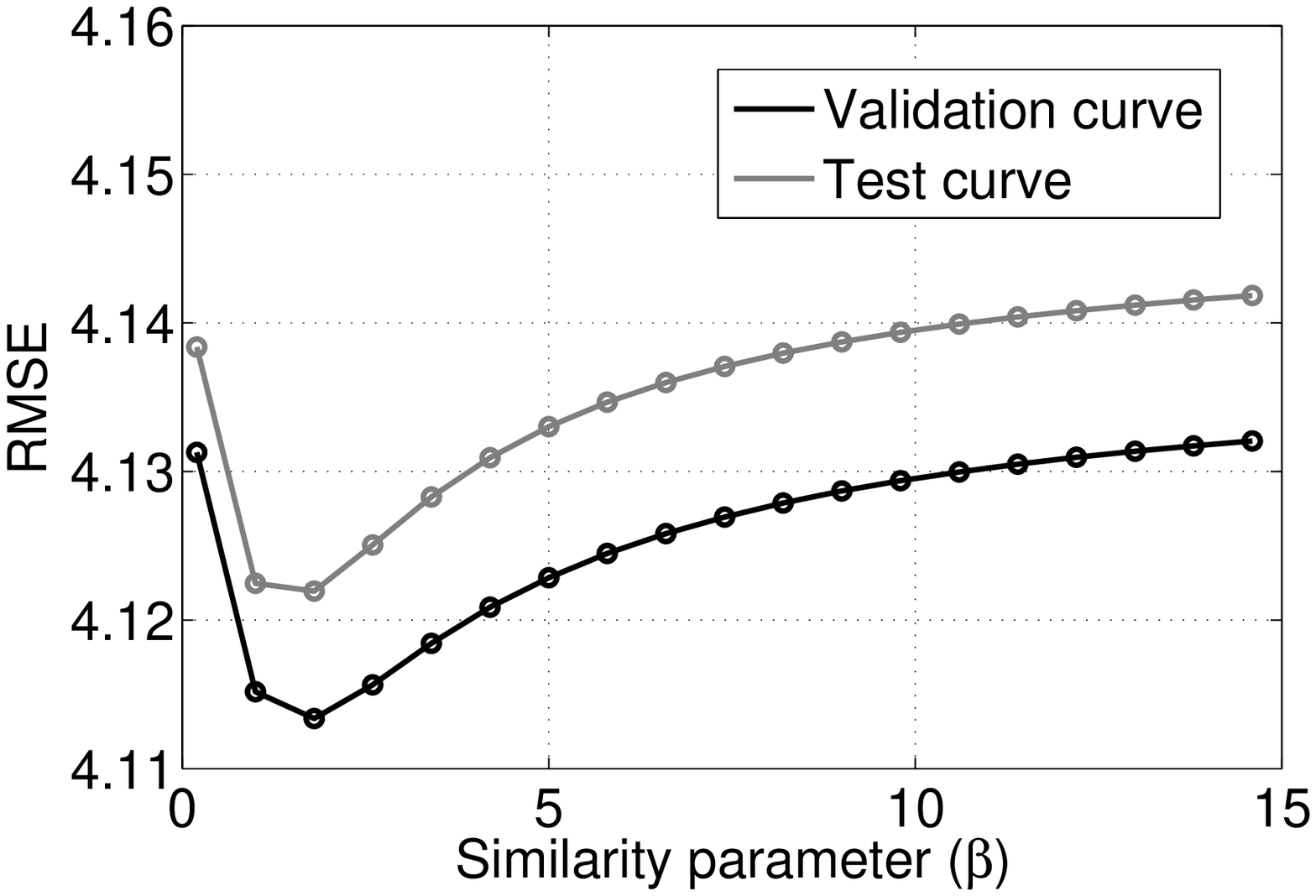}}
\subfloat[][]{\includegraphics[width=6.05cm]{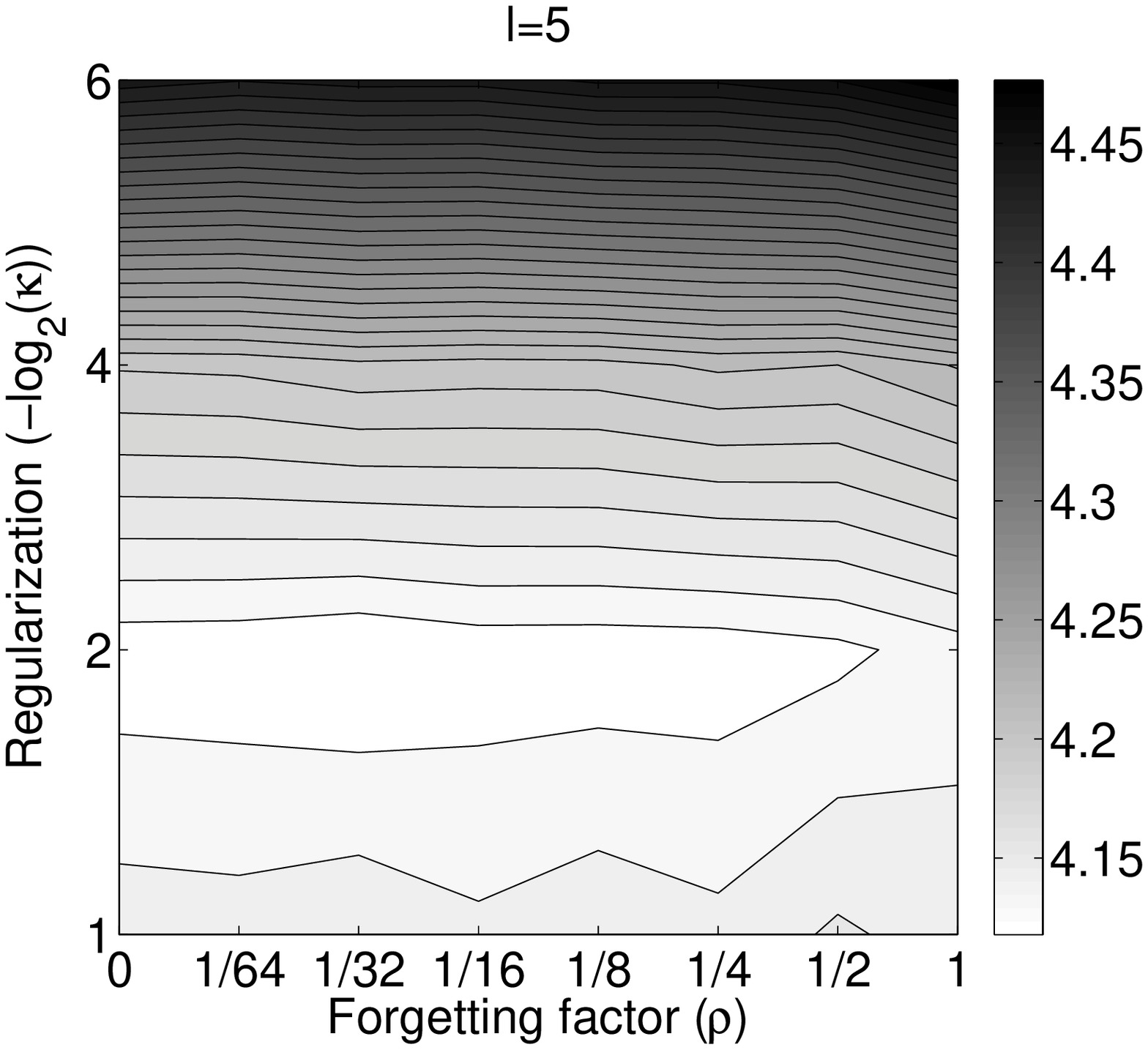}}
\subfloat[][]{\includegraphics[width=6.05cm]{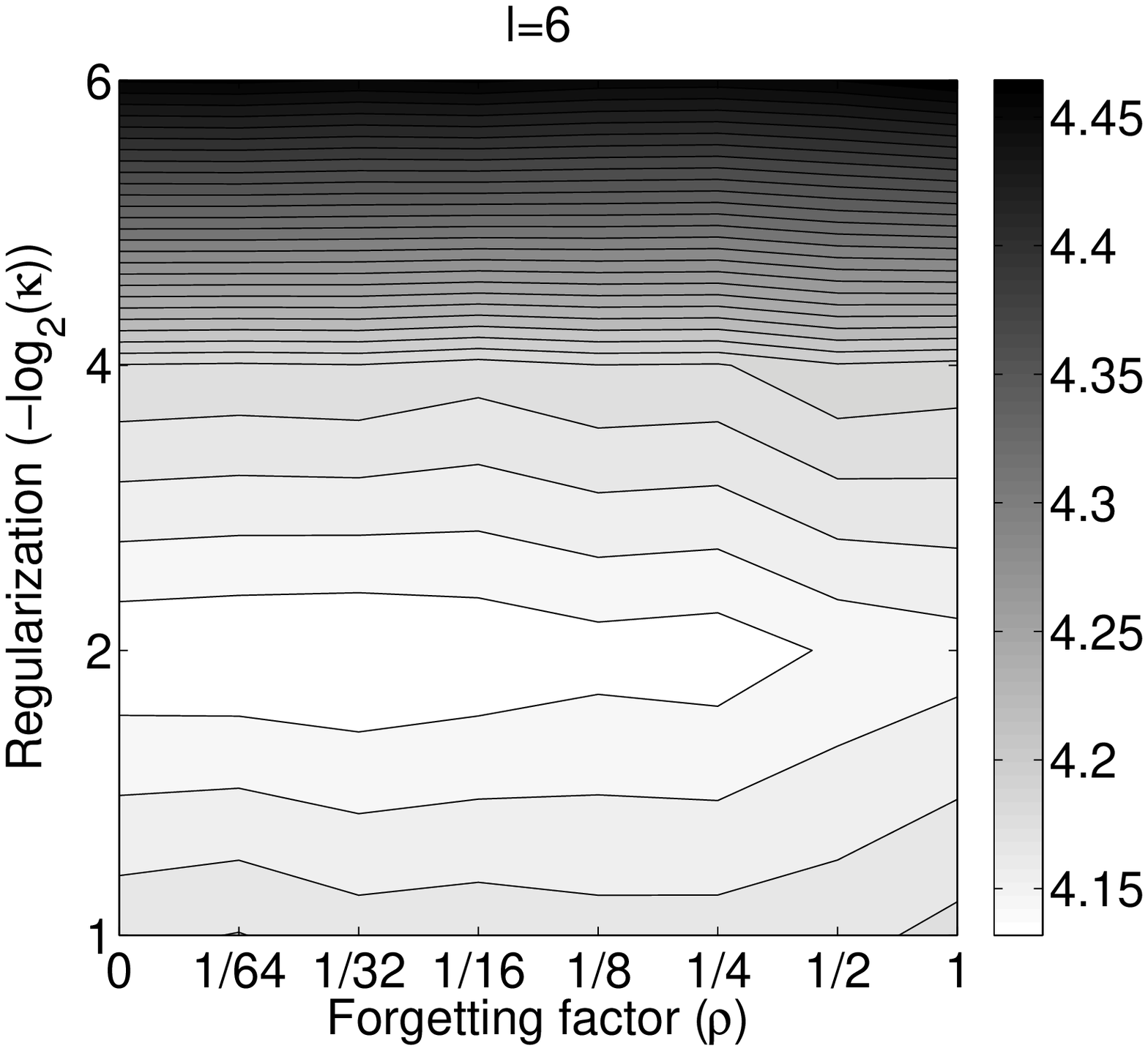}}
\caption[]{RMSE validation surfaces [(a)-(b), (e)-(f)] and test surfaces (c) as a function of forgetting factor ($\rho$) and regularization ($\kappa$).
(d): validation and test curve using the optimal number of levels $l=4$, regularization weight $\kappa=\frac{1}{2^{2}}$, forgetting
factor $\rho=0$, mini-bach size $R=8$, similarity parameter $\beta=1.8$. Group structure ($\G$): complete binary tree. Neighbor correction: $S_1^0$.
(a)-(c),(e)-(f): for fixed $(\kappa, \rho)$ parameter pair, the surfaces show the best RMSE values optimized in the $\beta$ similarity parameter.
(a): $l=3$. (b)-(c): $l=4$, on the same scale. (e): $l=5$. (f): $l=6$.}
\label{fig:bintree:validation surfaces,val-vs-test curve}
\end{figure*}

\begin{figure*}%
\centering%
\subfloat[][]{\includegraphics[width=6.05cm]{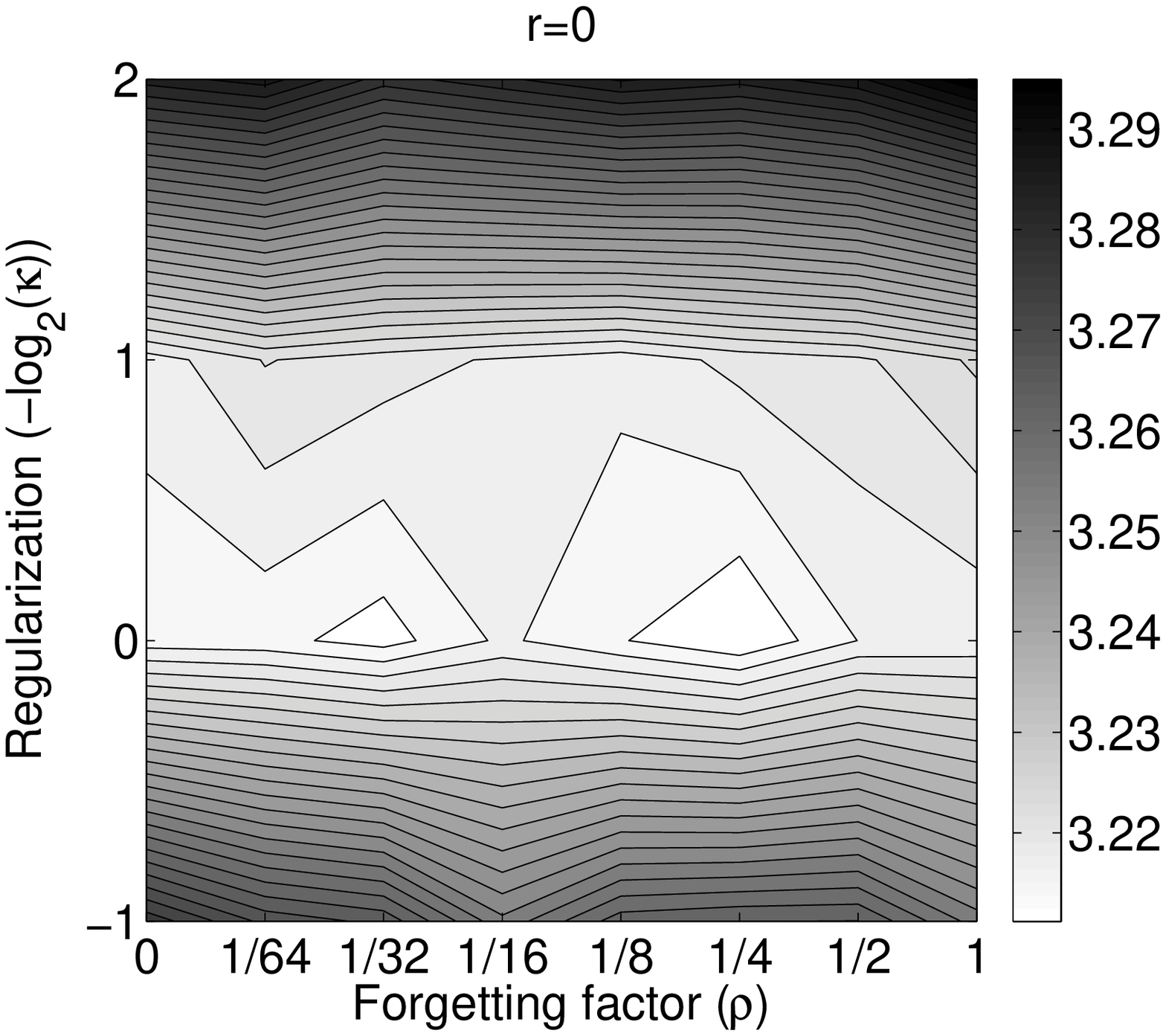}}
\subfloat[][]{\includegraphics[width=6.05cm]{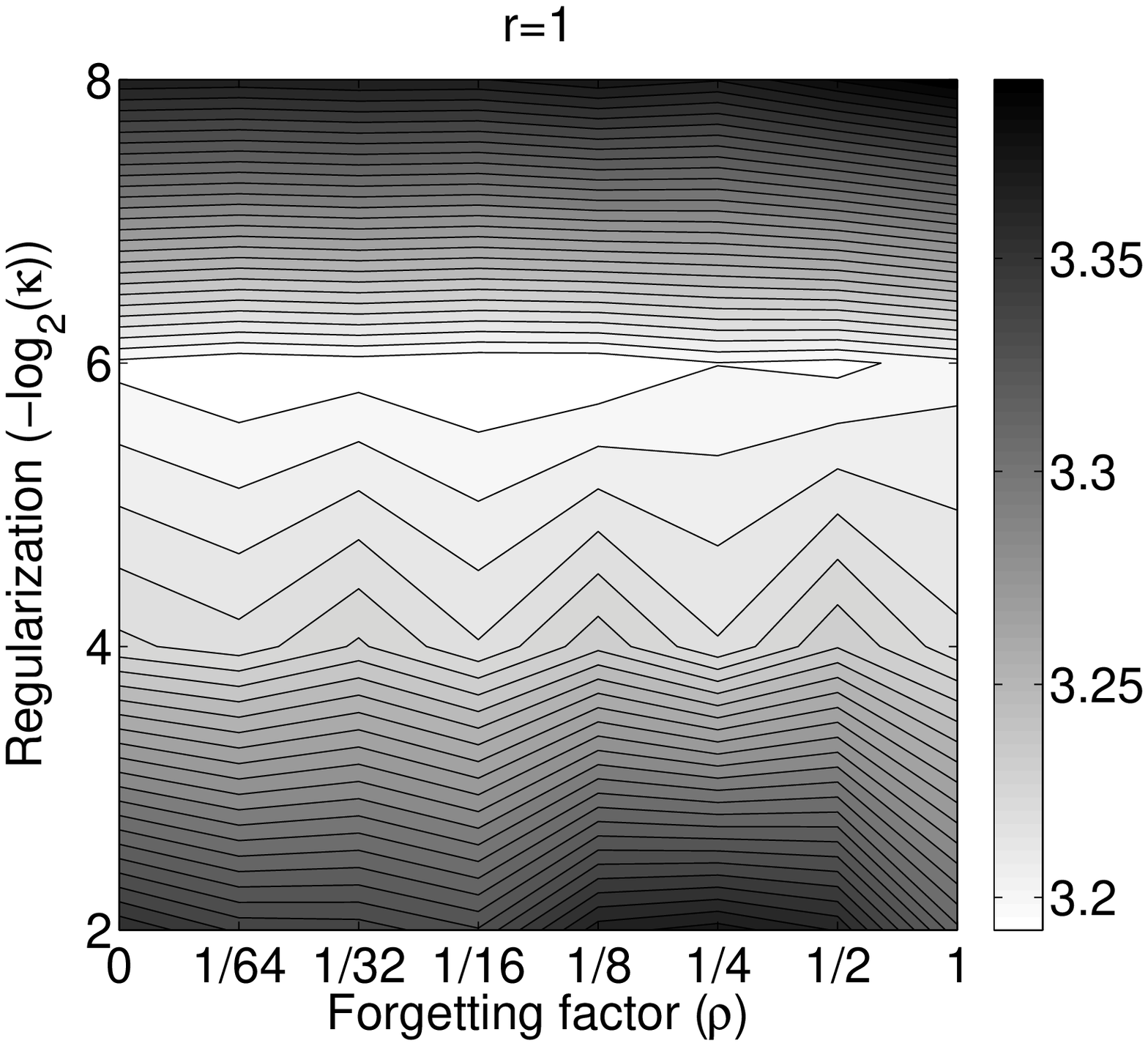}}
\subfloat[][]{\includegraphics[width=6.05cm]{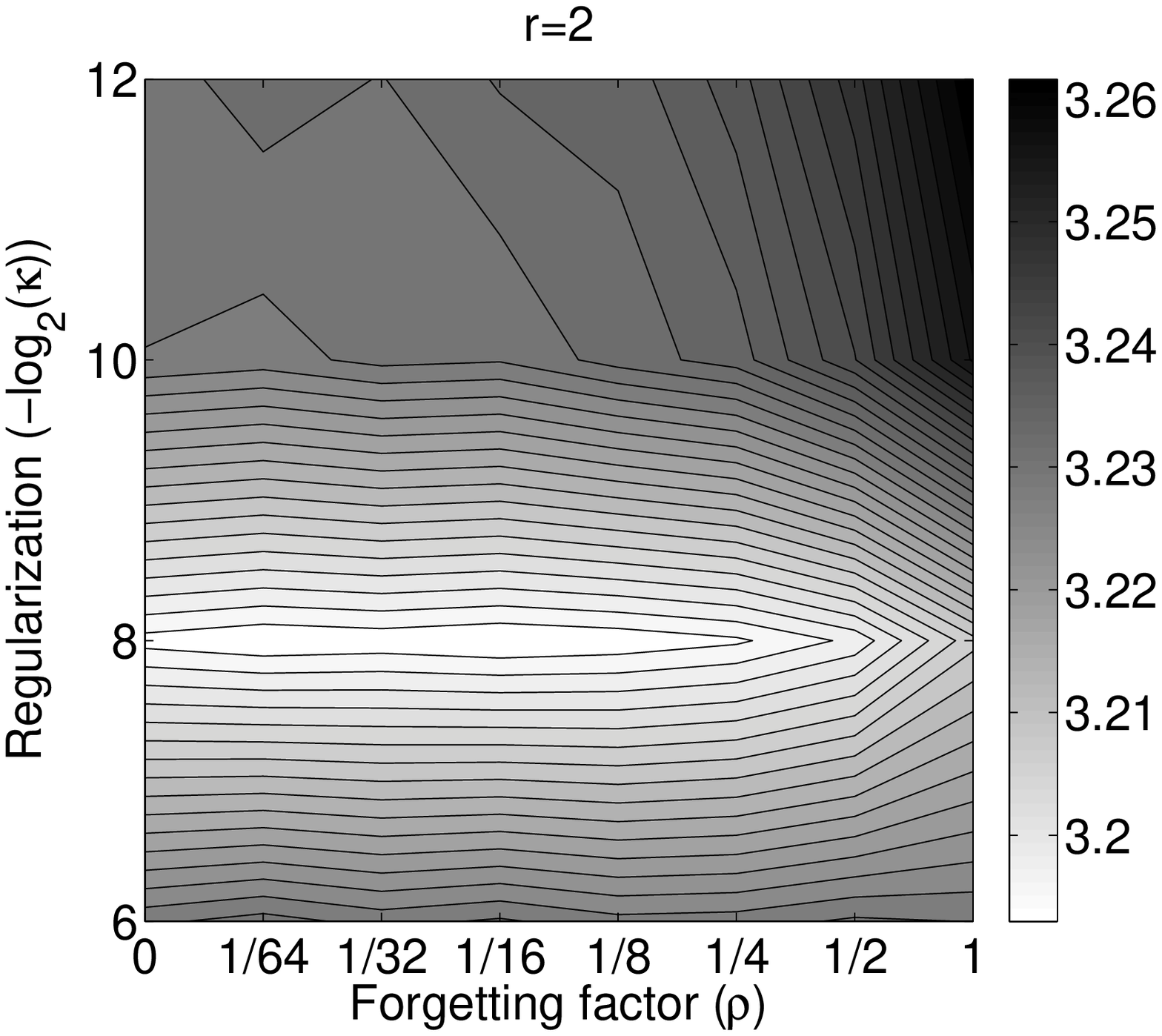}}\\
\subfloat[][]{\includegraphics[width=6.05cm]{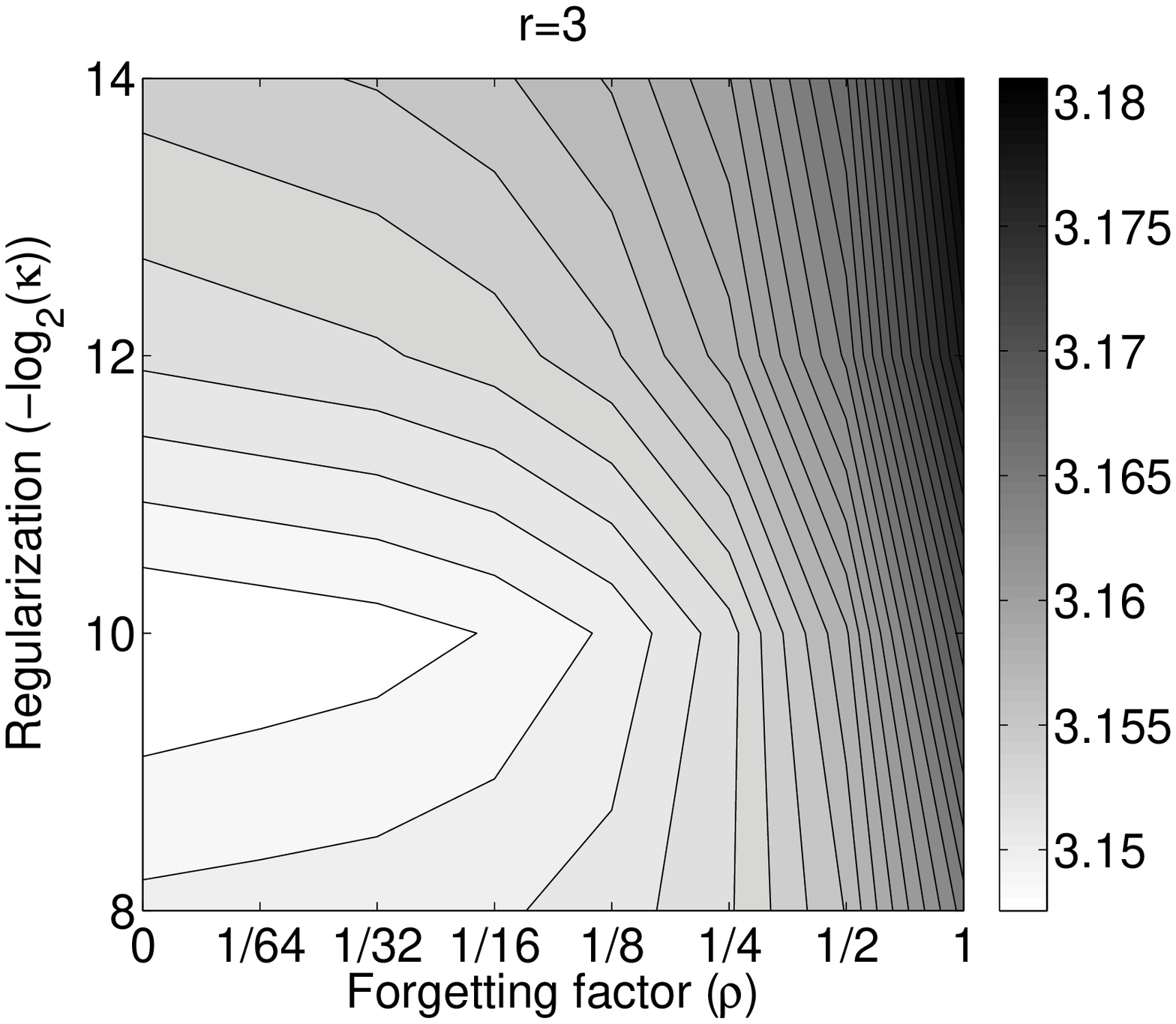}}
\subfloat[][]{\includegraphics[width=6.05cm]{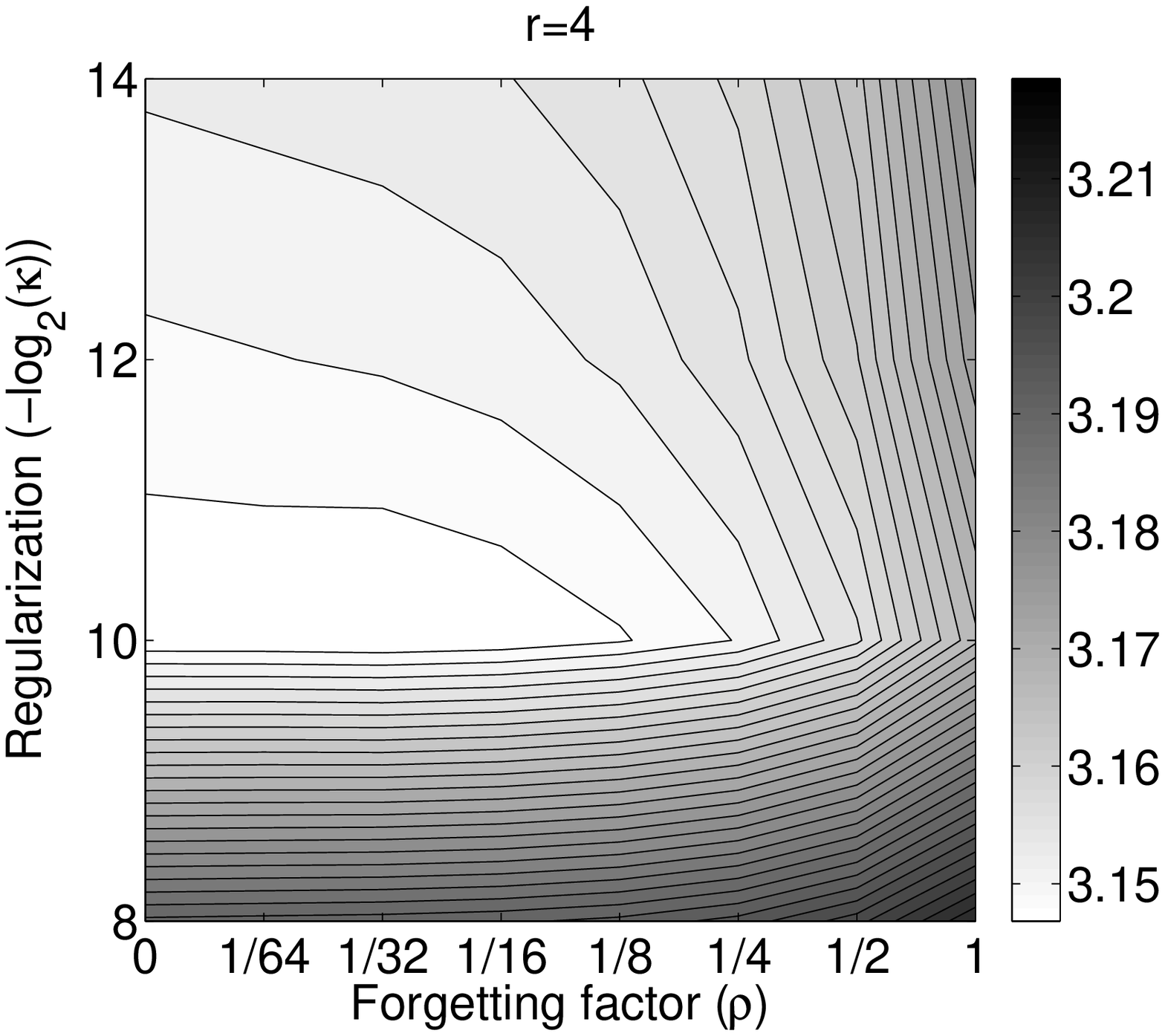}}
\subfloat[][]{\includegraphics[width=6.05cm]{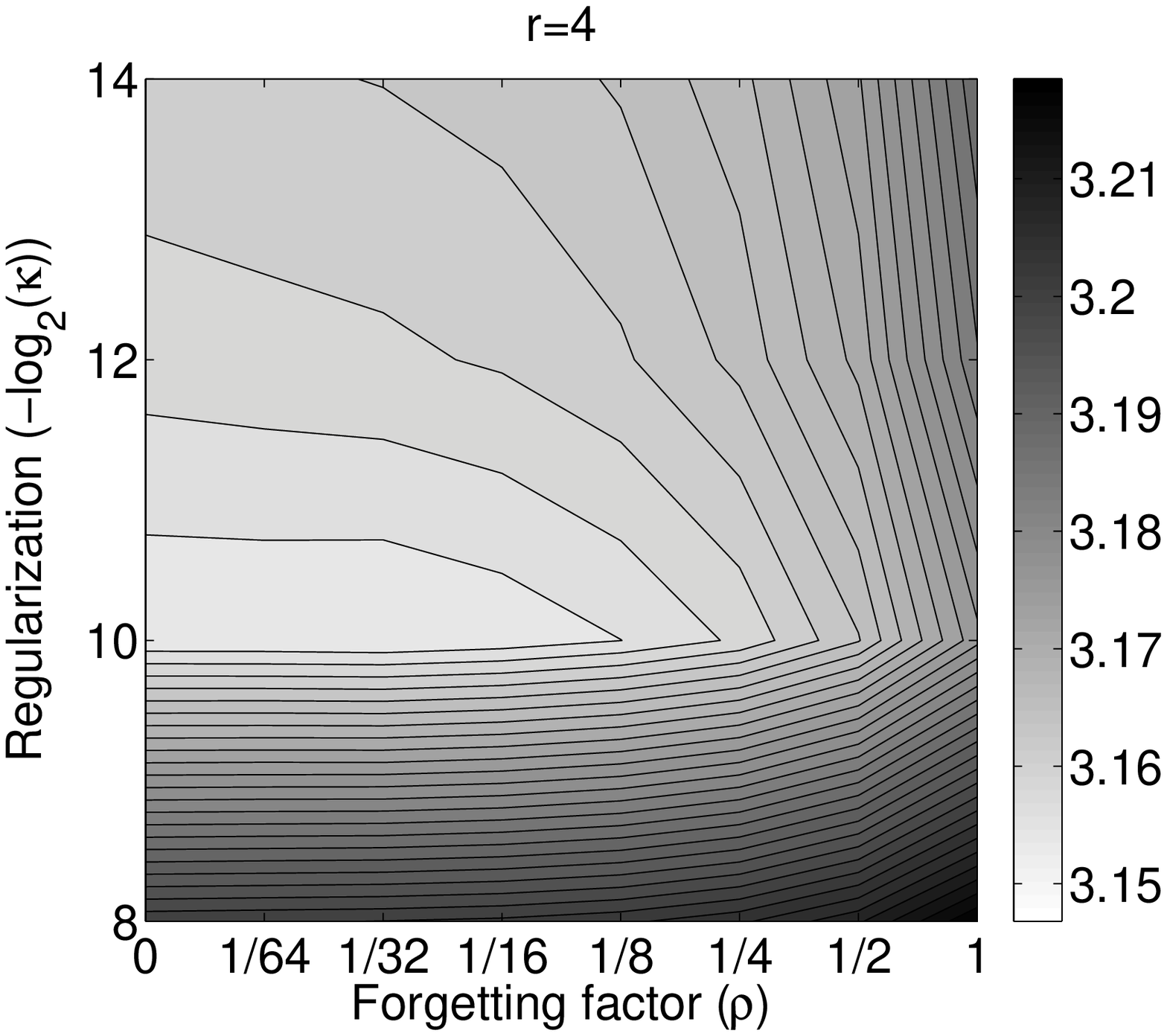}}
\caption[]{MAE validation surfaces [(a)-(e)] and test surfaces  (f) as a function of forgetting factor ($\rho$) and regularization ($\kappa$).
 For a fixed $(\kappa, \rho)$ parameter pair, the surfaces show the best MAE values optimized in the $\beta$ similarity parameter. The group structure ($\G$) is  toroid. The applied neighbor correction was $S_1^0$. (a): $r=0$ (no structure). (b): $r=1$. (c): $r=2$. (d): $r=3$. (e)-(f): $r=4$, on the same scale.}%
\label{fig:torus:MAE:validation surfaces}
\end{figure*}

\begin{figure}%
\centering%
\includegraphics[width=6.2cm]{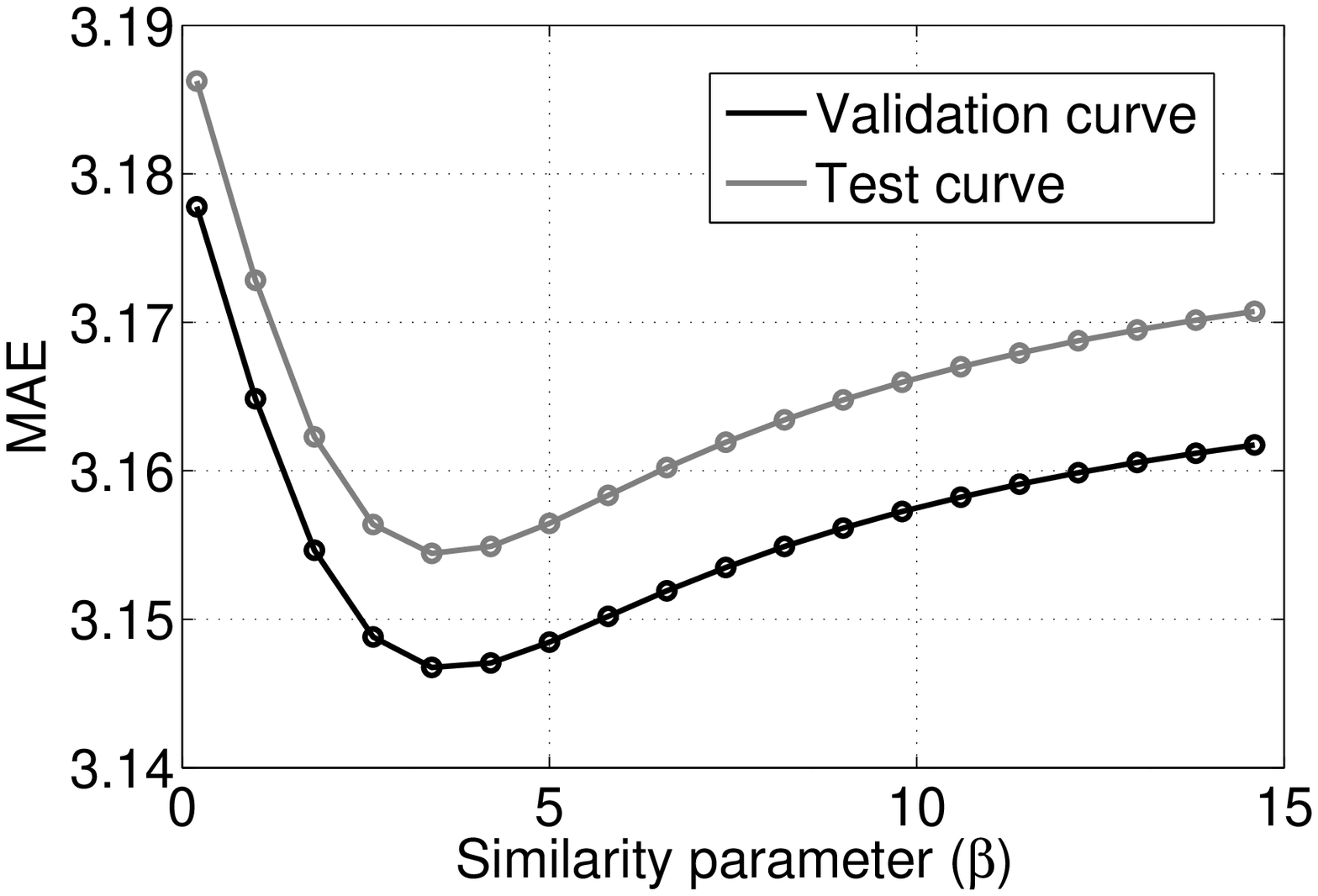}
\caption[]{MAE validation and test curves for toroid group structure using the optimal neighborhood size $r=4$, regularization weight $\kappa=\frac{1}{2^{10}}$, forgetting factor $\rho=\frac{1}{2^5}$, mini-batch size $R=8$, and similarity parameter $\beta=3.4$. The applied neighbor correction was $S_1^0$.}%
\label{fig:torus:MAE:validation vs test curve}
\end{figure}

When we used the \emph{MAE performance} measure, our results were similar to those of the RMSE.
We got the best results using toroid group structure, thus we present more details for this case.
\begin{itemize}
 \item With the usage of structured dictionaries we can get better results: the estimation errors were decreasing when we increased the neighbor size $r$ up to 4. (Table~\ref{tab:MAE,torus:perf:r}).
       The validation and test surfaces/curves are very similar, see Fig.~\ref{fig:torus:MAE:validation surfaces}(e)-(f), Fig.~\ref{fig:torus:MAE:validation vs test curve}.
 \item The quality of the estimation depends mostly on the $\kappa$ regularization parameter (Fig.~\ref{fig:torus:MAE:validation surfaces}(a)-(e)). The
	applied $\rho$ forgetting factor, $R$ mini-batch size and neighbor correction method can help in fine-tuning the results, see
      Fig.~\ref{fig:torus:MAE:validation surfaces}(a)-(e), Table~\ref{tab:MAE,torus:perf:r} and Table~\ref{tab:MAE,torus:perf:S}, respectively.
\item The smallest MAE we could achieve was $3.1544$, using $r=4$ neighbor size, $S_1^0$ neighbor correction and  $R=8$ mini-batch size.
      The baseline methods led to \cite{takacs08matrix} MAE = $3.1616$, \cite{takacs09scalable} MAE = $3.1606$ results. Our approach outperformed both of the state-of-the-art competitors. We also repeated this experiment $5$ more times with different randomly selected training, test, and validation sets,
      and our MAE results have never been worse than $3.155$. This demonstrates the efficiency of our approach.
\end{itemize}

\begin{table}
    \centering
    \caption{
    Performance (MAE) of the OSDL prediction using toroid group structure ($\G$) with different neighbor sizes $r$ ($r=0$: unstructured case). First-second row: mini-batch size $R=8$, third-fourth row: $R=16$.
        Odd rows: $S_1^0$, even rows: $S_2^0$ neighbor correction.  For fixed $R$, the best performance is highlighted with boldface typesetting.}
    \begin{tabular}{|c|c|c|c|c|c|c|}
    \hline
       &       &  $r=0$ & $r=1$ & $r=2$ & $r=3$ & $r=4$ \\
    \hline\hline
     $R=8$ & $S_1^0$ & $3.2225$ & $3.2019$ & $3.1989$ & $3.1563$ &  $\mathbf{3.1544}$\\
     \hline
           & $S_2^0$ & $3.2371$ & $3.2151$ & $3.2085$ &  $3.1584$&  $3.1571$\\
     \hline\hline
     $R=16$ & $S_1^0$ & $3.2220$ & $3.1988$ & $3.1982$ & $3.1576$ & $\mathbf{3.1546}$\\
    \hline
            & $S_2^0$ & $3.2382$ & $3.2147$ & $3.2101$ & $3.1594$ &  $3.1568$\\
    \hline
    \end{tabular}
    \label{tab:MAE,torus:perf:r}
\end{table}

\begin{table}
    \caption{Performance (MAE) of the OSDL prediction for different neighbor corrections
    using toroid group structure ($\G$). Columns: applied neighbor corrections. Rows: mini-batch size $R=8$ and $16$.
     The neighbor size was set to $r=4$. For fixed $R$, the best performance is highlighted with boldface typesetting.}
    \centering
    \begin{tabular}{|c|c|c|c|c|}
    \hline
          & $S_1$ & $S_2$ & $S_1^0$ & $S_2^0$\\
   \hline\hline
    $R=8$ & $3.1719$ & $3.1779$ & $\mathbf{3.1544}$ & $3.1571$\\
   \hline
    $R=16$ & $3.1726$ & $3.1778$ & $\mathbf{3.1546}$ & $3.1568$\\
    \hline
    \end{tabular}
    \label{tab:MAE,torus:perf:S}
\end{table}

\section{Conclusions}\label{sec:conclusions}
We have dealt with collaborative filtering (CF) based recommender systems and extended the application
domain of structured dictionaries to CF. We used online group-structured dictionary learning (OSDL) to solve the CF problem;
we casted the CF estimation task as an OSDL problem. We demonstrated the applicability of our novel approach on
joke recommendations. Our extensive numerical experiments show that structured dictionaries have several advantages
over the state-of-the-art CF methods: more precise estimation can be obtained, and smaller dimensional feature representation can be sufficient by applying
group structured dictionaries. Moreover, the estimation behaves robustly as a function of the OSDL parameters and the applied group structure.

\section*{Acknowledgments} The Project is supported by the European Union and co-financed by the
European Social Fund (grant agreements no.\ T\'AMOP 4.2.1/B-09/1/KMR-2010-0003 and KMOP-1.1.2-08/1-2008-0002). The
research was partly supported by the Department of Energy (grant number DESC0002607).

\bibliographystyle{IEEEtran}
%\bibliography{BIBs/IEEEabrv,BIBs/CF_via_OSDL_short}

\begin{thebibliography}{10}
\providecommand{\url}[1]{#1}
\csname url@samestyle\endcsname
\providecommand{\newblock}{\relax}
\providecommand{\bibinfo}[2]{#2}
\providecommand{\BIBentrySTDinterwordspacing}{\spaceskip=0pt\relax}
\providecommand{\BIBentryALTinterwordstretchfactor}{4}
\providecommand{\BIBentryALTinterwordspacing}{\spaceskip=\fontdimen2\font plus
\BIBentryALTinterwordstretchfactor\fontdimen3\font minus
  \fontdimen4\font\relax}
\providecommand{\BIBforeignlanguage}[2]{{%
\expandafter\ifx\csname l@#1\endcsname\relax
\typeout{** WARNING: IEEEtran.bst: No hyphenation pattern has been}%
\typeout{** loaded for the language `#1'. Using the pattern for}%
\typeout{** the default language instead.}%
\else
\language=\csname l@#1\endcsname
\fi
#2}}
\providecommand{\BIBdecl}{\relax}
\BIBdecl

\bibitem{ricci11recommender}
F.~Ricci, L.~Rokach, B.~Shapira, and P.~Kantor, \emph{Recommender Systems
  Handbook}.\hskip 1em plus 0.5em minus 0.4em\relax Springer, 2011.

\bibitem{takacs09scalable}
G.~Tak{\'a}cs, I.~Pil{\'a}szy, B.~N{\'e}meth, and D.~Tikk, ``Scalable
  collaborative filtering approaches for large recommender systems,'' \emph{J.
  Mach. Learn. Res.}, vol.~10, pp. 623--656, 2009.

\bibitem{yaghoobi09dictionary}
M.~Yaghoobi, T.~Blumensath, and M.~Davies, ``Dictionary learning for sparse
  approximations with the majorization method,'' \emph{{IEEE} Trans. Signal
  Process.}, vol.~57, no.~6, pp. 2178--2191, 2009.

\bibitem{witten09penalized}
D.~M. Witten, R.~Tibshirani, and T.~Hastie, ``A penalized matrix decomposition,
  with applications to sparse principal components and canonical correlation
  analysis,'' \emph{Biostatistics}, vol.~10, no.~3, pp. 515--534, 2009.

\bibitem{zou06sparse}
H.~Zou, T.~Hastie, and R.~Tibshirani, ``Sparse principal component analysis,''
  \emph{J. Comput. Graph. Stat.}, vol.~15, no.~2, pp. 265--286, 2006.

\bibitem{hyvarinen01independent}
A.~Hyv{\"a}rinen, J.~Karhunen, and E.~Oja, \emph{Independent Component
  Analysis}.\hskip 1em plus 0.5em minus 0.4em\relax John Wiley \& Sons, 2001.

\bibitem{cardoso98multidimensional}
J.~Cardoso, ``Multidimensional independent component analysis,'' in
  \emph{ICASSP 1998}, pp. 1941--1944.

\bibitem{lee00algorithms}
D.~D. Lee and H.~S. Seung, ``Algorithms for non-negative matrix
  factorization,'' in \emph{NIPS 2000}, pp. 556--562.

\bibitem{tropp10computational}
J.~A. Tropp and S.~J. Wright, ``Computational methods for sparse solution of
  linear inverse problems,'' \emph{Proc. of the IEEE special issue on
  Applications of sparse representation and compressive sensing}, vol.~98,
  no.~6, pp. 948--958, 2010.

\bibitem{natarajan95sparse}
B.~K. Natarajan, ``Sparse approximate solutions to linear systems,'' \emph{SIAM
  J. Comput.}, vol.~24, no.~2, pp. 227--234, 1995.

\bibitem{tibshirani96regression}
R.~Tibshirani, ``Regression shrinkage and selection via the {Lasso},'' \emph{J.
  Roy. Stat. Soc. B. Met.}, vol.~58, no.~1, pp. 267--288, 1996.

\bibitem{huang10benefit}
J.~Huang and T.~Zhang, ``The benefit of group sparsity,'' \emph{Ann. Stat.},
  vol.~38, no.~4, pp. 1978--2004, 2010.

\bibitem{yuan06model}
M.~Yuan and Y.~Lin, ``Model selection and estimation in regression with grouped
  variables,'' \emph{J. Roy. Stat. Soc. B Met.}, vol.~68, no.~1, pp. 49--67,
  2006.

\bibitem{baraniuk10model}
R.~G. Baraniuk, V.~Cevher, M.~F. Duarte, and C.~Hegde, ``Model-based
  compressive sensing,'' \emph{{IEEE} Trans. Inf. Theory}, vol.~56, pp. 1982 --
  2001, 2010.

\bibitem{zhang10automatic}
S.~Zhang, J.~Huang, Y.~Huang, Y.~Yu, H.~Li, and D.~Metaxas, ``Automatic image
  annotation using group sparsity,'' in \emph{CVPR 2010}, pp. 3312--3319.

\bibitem{zhao09composite}
P.~Zhao, G.~Rocha, and B.~Yu, ``The composite absolute penalties family for
  grouped and hierarchical variable selection,'' \emph{Ann. Stat.}, vol.~37,
  no.~6A, pp. 3468--3497, 2009.

\bibitem{jacob09group}
L.~Jacob, G.~Obozinski, and J.-P. Vert, ``Group {Lasso} with overlap and graph
  {Lasso},'' in \emph{ICML 2009}, pp. 433--440.

\bibitem{kim10tree}
S.~Kim and E.~P. Xing, ``Tree-guided group {Lasso} for multi-task regression
  with structured sparsity,'' in \emph{ICML 2010}, pp. 543--550.

\bibitem{rapaport08classification}
F.~Rapaport, E.~Barillot, and J.-P. Vert, ``Classification of array{CGH} data
  using fused {SVM},'' \emph{Bioinformatics}, vol.~24, pp. i375--i382, 2008.

\bibitem{obozinski10joint}
G.~Obozinski, B.~Taskar, and M.~I. Jordan, ``Joint covariate selection and
  joint subspace selection for multiple classification problems,'' \emph{Stat.
  Comput.}, vol.~20, pp. 231--252, 2010.

\bibitem{kim10scalable}
D.~Kim, S.~Sra, and I.~S. Dhillon, ``A scalable trust-region algorithm with
  application to mixed-norm regression,'' in \emph{ICML 2010}, pp. 519--526.

\bibitem{rakotomamonjy11review}
A.~Rakotomamonjy, ``Review: Surveying and comparing simultaneous sparse
  approximation (or group-lasso) algorithms,'' \emph{Signal Process.}, vol.~91,
  no.~7, pp. 1505--1526, 2011.

\bibitem{szafranski10composite}
M.~Szafranski, Y.~Grandvalet, and A.~Rakotomamonjy, ``Composite kernel
  learning,'' \emph{Mach. Learn.}, vol.~79, pp. 73--103, 2010.

\bibitem{aflalo11variable}
J.~Aflalo, A.~Ben-Tal, C.~Bhattacharyya, J.~S. Nath, and S.~Raman, ``Variable
  sparsity kernel learning,'' \emph{J. Mach. Learn. Res.}, vol.~12, pp.
  565--592, 2011.

\bibitem{elhamifar11robust}
E.~Elhamifar and R.~Vidal, ``Robust classification using structured sparse
  representation,'' in \emph{CVPR 2011}, pp. 1873 -- 1879.

\bibitem{schmidt10convex}
M.~Schmidt and K.~Murphy, ``Convex structure learning in log-linear models:
  Beyond pairwise potentials,'' \emph{AISTATS, J. Mach. Learn. Res.:W\&CP},
  vol.~9, pp. 709--716, 2010.

\bibitem{jalali11learning}
A.~Jalali, P.~Ravikumar, V.~Vasuki, and S.~Sanghavi, ``On learning discrete
  graphical models using group-sparse regularization,'' \emph{AISTATS,
  JMLR:W\&CP}, vol.~15, 2011.

\bibitem{bach11optimization}
F.~Bach, R.~Jenatton, J.~Marial, and G.~Obozinski, \emph{Optimization for
  Machine Learning}.\hskip 1em plus 0.5em minus 0.4em\relax MIT Press, 2011,
  ch. Convex optimization with sparsity-inducing norms.

\bibitem{jenatton11proximal}
R.~Jenatton, J.~Mairal, G.~Obozinski, and F.~Bach, ``Proximal methods for
  hierarchical sparse coding,'' \emph{J. Mach. Learn. Res.}, vol.~12, pp.
  2297--2334, 2011.

\bibitem{jenatton10structured}
R.~Jenatton, G.~Obozinski, and F.~Bach, ``Structured sparse principal component
  analysis,'' \emph{AISTATS, J. Mach. Learn. Res.:W\&CP}, vol.~9, pp. 366--373,
  2010.

\bibitem{mairal11convex}
J.~Mairal, R.~Jenatton, G.~Obozinski, and F.~Bach, ``Convex and network flow
  optimization for structured sparsity,'' \emph{J. Mach. Learn. Res.}, vol.~12,
  pp. 2681--2720, 2011.

\bibitem{rosenblum10dictionary}
K.~Rosenblum, L.~Zelnik-Manor, and Y.~Eldar, ``Dictionary optimization for
  block-sparse representations,'' in \emph{AAAI Fall 2010 Symposium on Manifold
  Learning}.

\bibitem{koray09learning}
K.~Kavukcuoglu, M.~Ranzato, R.~Fergus, and Y.~LeCun, ``Learning invariant
  features through topographic filter maps,'' in \emph{CVPR 2009}, pp.
  1605--1612.

\bibitem{silva11blind}
J.~Silva, M.~Chen, Y.~C. Eldar, G.~Sapiro, and L.~Carin, ``Blind compressed
  sensing over a structured union of subspaces,'' Tech. Rep., 2011,
  \url{http://arxiv.org/abs/1103.2469}.

\bibitem{bottou05on-line}
L.~Bottou and Y.~L. Cun, ``On-line learning for very large data sets,''
  \emph{Appl. Stoch. Model. Bus. - Stat. Learn.}, vol.~21, no.~2, pp. 137--151,
  2005.

\bibitem{szabo11online}
Z.~Szab{\'o}, B.~P{\'o}czos, and A.~L{\H{o}}rincz, ``Online group-structured
  dictionary learning,'' in \emph{CVPR 2011}, pp. 2865--2872.

\bibitem{bertsekas99nonlinear}
D.~P. Bertsekas, \emph{Nonlinear Programming}.\hskip 1em plus 0.5em minus
  0.4em\relax Athena Scientific Belmont, 1999.

\bibitem{goldber01eigentaste}
K.~Goldberg, T.~Roeder, D.~Gupta, and C.~Perkins, ``Eigentaste: A constant time
  collaborative filtering algorithm,'' \emph{Inform. Retrieval}, vol.~4, pp.
  133--151, 2001.

\bibitem{takacs08matrix}
G.~Tak{\'a}cs, I.~Pil{\'a}szy, B.~N{\'e}meth, and D.~Tikk, ``Matrix
  factorization and neighbor based algorithms for the {Netflix} prize
  problem,'' in \emph{RecSys 2008}, pp. 267--274.

\end{thebibliography}

\end{document}